\documentclass[11pt]{article}
\usepackage{graphicx, color}
\usepackage{amsmath,amsfonts,amssymb}
\usepackage{graphicx}

\def\ds{\displaystyle}
\def\forall{\hbox{for all}~}
\def\L{{\bf L}}

\def\argmax{\hbox{arg}\!\max}
\def\bfv{{\bf v}}
\def\bfx{{\bf x}}
\def\bfw{{\bf w}}

\def\bfy{{\bf y}}
\def\bfn{{\bf n}}
\def\bfe{{\bf e}}

\def\bfS{{\bf S}}
\def\ve{\varepsilon}

\def\E{{\cal E}}
\def\I{{\cal I}}
\def\J{{\cal J}}
\def\S{{\cal S}}

\def\R{{\mathbb R}}

\def\implies{\Longrightarrow}
\def\vp{\varphi}

\def\vs{\vskip 2em}

\def\v{\vskip 1em}
\def\O{{\cal O}}

\def\begi{\begin{itemize}}
\def\endi{\end{itemize}}

\def\B{{\cal B}}

\def\ov{\overline}
\def\Tilde{\widetilde}

\def\bega{\begin{array}}
\def\enda{\end{array}}
\def\meas{\hbox{meas}}

\def\bel{\begin{equation}\label}
\def\eeq{\end{equation}}

\def\sqr#1#2{\vbox{\hrule height .#2pt
\hbox{\vrule width .#2pt height #1pt \kern #1pt
\vrule width .#2pt}\hrule height .#2pt }}
\def\square{\sqr74}
\def\endproof{\hphantom{MM}\hfill\llap{$\square$}\goodbreak}

\newtheorem{theorem}{Theorem}[section]
\newtheorem{lemma}[theorem]{Lemma}

\newtheorem{corollary}[theorem]{Corollary}
\newtheorem{definition}[theorem]{Definition}
\newtheorem{remark}[theorem]{Remark}

\begin{document}

\title{\bf On a Shape Optimization Problem for  Tree Branches}
\vs

\author{Alberto Bressan$^{(*)}$ and  Sondre T.~Galtung$^{(**)}$\\
\,
\\
$^{(*)}$ Department of Mathematics, Penn State University \\
University Park, Pa.~16802, USA.\\
\,
\\
$^{(**)}$ Department of Mathematical Sciences, \\
NTNU -- Norwegian University of Science and Technology, \\ NO-7491 Trondheim, Norway.
\\
\,
\\
e-mails: axb62@psu.edu, sondre.galtung@ntnu.no.
}
\maketitle
\begin{abstract} This paper is concerned with a shape optimization 
problem, where
the functional to be maximized describes the total sunlight collected by 
a distribution of tree leaves, minus the cost for transporting water and nutrient from the 
base of the trunk to all the leaves.  In the  case of
2 space dimensions,  the solution is proved to be unique, and explicitly determined. 
\end{abstract}


\vs
\section{Introduction}
\label{s:0}
\setcounter{equation}{0}
In the recent papers \cite{BPS, BSu} two functionals were introduced, 
measuring the amount of light collected by the leaves, 
and the amount of water and nutrients collected by the roots of a tree.
In connection with a ramified transportation cost \cite{BCM, MMS, X03},
these lead to various optimization problems for tree shapes.

Quite often, optimal solutions to problems involving a
ramified transportation cost exhibit a fractal structure \cite{BCM1, BraS,
BW, DS,  MoS, PSX, S1}.  In the present note we analyze in more detail the 
optimization 
problem  for tree branches proposed in \cite{BPS}, in the
2-dimensional case.   In this simple setting,  the unique 
solution can be explicitly determined.  Instead of being fractal, 
its shape reminds of a solar panel.    

\begin{figure}[ht]
	\centerline{\hbox{\includegraphics[width=8cm]{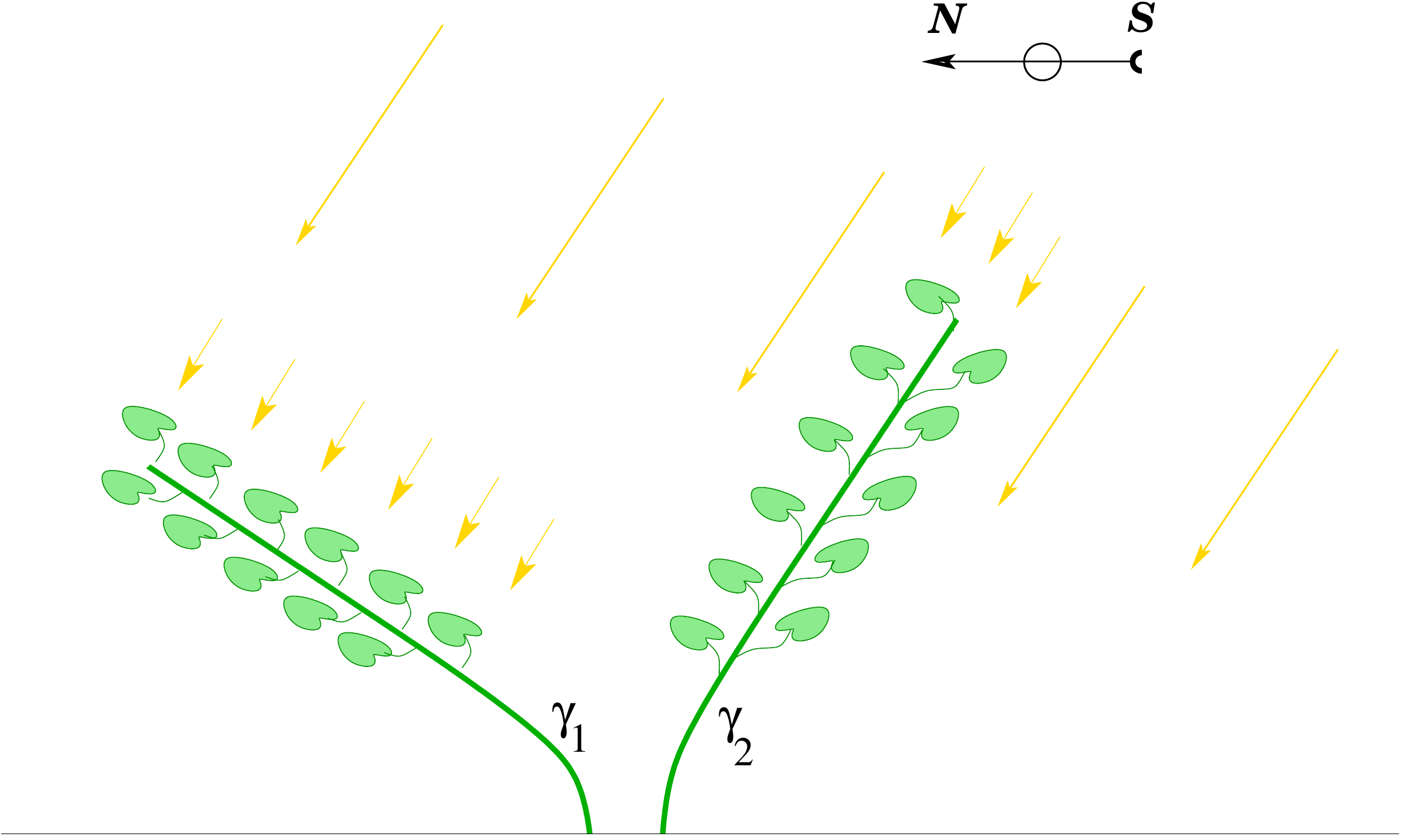}}}
	\caption{\small A stem $\gamma_1$ perpendicular to the sun rays is optimally shaped to collect the most
	light. For the stem $\gamma_2$ bending toward the light source, the upper leaves put the lower ones in shade.  
	} \label{f:pg45}
\end{figure}

The present analysis was partially motivated 
by the goal of understanding phototropism, i.e., the tendency of plant stems
to bend toward the source of light.   Our results indicate  that this behavior cannot be explained 
purely in terms of maximizing the amount of light collected by the leaves (Fig.~\ref{f:pg45}).
Apparently, other factors must have played a role in the evolution of this trait, 
such as the competition among different plants. See \cite{BGRR} for 
some results in this direction.

The remainder of this paper is organized as follows. In Section~2 we review the two functionals
defining the shape optimization problem, and state the main results.   Proofs are then worked out in 
Sections 3 to 5.

\v

\section{Statement of the main results}
\label{s:1}
\setcounter{equation}{0}

We begin by reviewing the two functionals considered in \cite{BPS, BSu}.

\subsection{A sunlight functional}


Let $\mu$ be a positive, bounded Radon measure on 
$\R^d_+\doteq\{(x_1,x_2,\ldots, x_d)\,;~x_d\geq 0\}$.   
Thinking of $\mu$ as the density of leaves on a tree, 
we seek a functional $\S(\mu)$ describing 
the total amount of sunlight absorbed by the leaves.
Fix a unit vector
$$\bfn~\in~S^{d-1}~\doteq~\{ x\in \R^d\,;~~|x|=1\},$$
and assume that all  
light rays come parallel to $\bfn$.  
Call $E_\bfn^\perp$ the $(d-1)$-dimensional subspace perpendicular to $\bfn$ 
and let $\pi_\bfn:\R^d\mapsto E_\bfn^\perp$ be the perpendicular projection.   Each point $\bfx\in \R^d$ can thus be expressed 
uniquely as
\bel{perp}
\bfx~=~\bfy + s\bfn\eeq
with $\bfy\in  E_\bfn^\perp$ and $s\in\R$.

On the perpendicular subspace $E_\bfn^\perp$ consider the projected measure $\mu^\bfn$, defined by setting
\bel{mupro}\mu^\bfn(A)~=~\mu\Big(\bigl\{ x\in\R^d\,;~~\pi_\bfn(x)\in A\bigr\}\Big).\eeq
Call $\Phi^\bfn$ the density of the absolutely continuous part of $\mu^\bfn$
w.r.t.~the $(d-1)$-dimensional Lebesgue measure on $E_\bfn^\perp$.
\v
\begin{definition}
The total amount of sunlight from the direction $\bfn$ captured by a measure 
$\mu$ on $\R^d$ is defined as
\bel{SSn}
\S^\bfn(\mu)~\doteq~
\int_{E_\bfn^\perp}\Big(1- \exp\bigl\{ - \Phi^\bfn(y)\bigr\}\Big)
\, dy\,.\eeq
More generally,
given an integrable function $\eta\in \L^1(S^{d-1})$, 
the total sunshine absorbed by $\mu$ from all directions
is defined as
\bel{SS2}
\S^\eta(\mu)~\doteq~
\int_{S^{d-1}}\left(\int_{E_\bfn^\perp}\Big(1- \exp\bigl\{ - \Phi^\bfn(y)\bigr\}\Big)
\, dy\right) \eta(\bfn)\,d\bfn\,.\eeq
\end{definition}
\v
In the formula (\ref{SS2}),    $\eta(\bfn)$ accounts for  the intensity of light
coming from  the direction $\bfn$.
\begin{remark}\label{r:22} {\rm According to the above definition, the amount of sunlight
$\S^\bfn(\mu)$
captured by the measure 
$\mu$ only depends on its projection $\mu^\bfn$ on the subspace perpendicular to 
$\bfn$.   In particular, if a second measure $\Tilde\mu$ is obtained from $\mu$
by shifting some of the mass in a direction parallel to $\bfn$, then 
$\S(\Tilde\mu) = \S(\mu)$.}
\end{remark}

\v

\subsection{Optimal irrigation patterns}
 
Consider a positive Radon measure $\mu$ on $\R^d$ with total mass
$M=\mu(\R^d)$, and let 
$\Theta=[0,M]$.
We think of  $\xi\in \Theta$ as a Lagrangian variable, labeling a water particle.

\begin{definition} A measurable map 
\bel{iplan}
\chi:\Theta\times \R_+~\mapsto~ \R^d\eeq
is called an {\bf admissible irrigation plan}
if 
\begi
\item[(i)] For every $\xi\in \Theta$, the map
$t\mapsto \chi(\xi,t)$ is Lipschitz continuous. 
More precisely, for each $\xi$ there exists a stopping time $T(\xi)$ such that, calling 
$$\dot \chi(\xi,t)~=~{\partial\over\partial t} ~\chi(\xi,t)$$
the partial derivative w.r.t.~time, one has
\bel{stime}\bigl|\dot \chi(\xi,t)\bigr|~=~\left\{ \bega{rl} 1\qquad &\quad\hbox{for a.e.}
~t\in \bigl[0, T(\xi)\bigr],\\[3mm]
0\qquad &\quad\hbox{for}
~t> T(\xi).\enda\right.\eeq
\item[(ii)] At time $t=0$ all particles are at the origin:
$\chi(\xi,0)={\bf 0}$ for all $\xi\in\Theta$.
\item[(iii)] The push-forward of the Lebesgue measure on $[0,M]$ through the map $\xi\mapsto 
\chi(\xi,T(\xi))$ coincides with the measure $\mu$.
In other words, for every open set $A\subset\R^d$ there holds
\bel{chi1}\mu(A)~=~\hbox{\rm meas}\Big( \{ \xi\in \Theta\,;~~\chi(\xi,T(\xi))\in A\bigr\}\Big).\eeq
\endi
\end{definition}

One may think of  $\chi(\xi,t)$ as the 
position of the water particle $\xi$ at time $t$.

To define the corresponding transportation cost, we first compute 
how many particles travel through a point $x\in\R^d$.  
This is described by
\bel{chi}|x|_\chi~\doteq~\meas\Big(\bigl\{\xi\in \Theta\,;~~\chi(\xi,t)= x~~~\hbox{for some}~~t\geq 0\bigr\}\Big).\eeq
We think of $|x|_\chi$ as the {\it total flux going through the
point $x$}.    Following \cite{G, MMS}, we consider
\v
\begin{definition}
{\bf (irrigation cost).} 
For a given $\alpha\in [0,1]$,
the total cost of the irrigation plan $\chi$ is
\bel{TCg}
\E^\alpha(\chi)~\doteq~\int_\Theta\left(\int_0^{T(\xi)} \bigl|\chi(\xi,t)
\bigr|_\chi^{\alpha-1} \, dt\right)
d\xi.\eeq
The  {\bf $\alpha$-irrigation cost} of a measure $\mu$
is defined as 
\bel{Idef}\I^\alpha(\mu)~\doteq~\inf_\chi \E^\alpha(\chi),\eeq
where the infimum is taken over all admissible irrigation plans for the measure $\mu$.
\end{definition}

\begin{remark} {\rm Sometimes it is convenient to consider more general
irrigation plans where, in place of (\ref{stime}),  for  a.e.~$t\in [0,T(\xi)]$ 
the speed satisfies $|\dot\chi(\xi,t)|\leq 1$.
In this case, the cost (\ref{TCg}) is replaced by
\bel{TC2}
\E^\alpha(\chi)~\doteq~\int_\Theta\left(\int_0^{T(\xi)} \bigl|\chi(\xi,t)
\bigr|_\chi^{\alpha-1} \,|\dot\chi(\xi,t)|\, dt\right)
d\xi.\eeq
Of course, one can always re-parameterize each trajectory 
$t\mapsto \chi(\xi,t)$ by arc-length, so that (\ref{stime}) holds.
This does not affect the cost (\ref{TC2}).
}
\end{remark}

\v
\begin{remark} {\rm In the case $\alpha=1$, the expression (\ref{TCg}) reduces to
$$
\E^\alpha(\chi)~\doteq~\int_\Theta\left(\int_{\R_+} |\dot \chi_t(\xi,t)|\, dt\right)
d\xi~=~\int_\Theta[\hbox{total length of the path} ~\chi(\xi,\cdot)]\, d\xi\,.$$
Of course, this length is minimal if every path $\chi(\cdot,\xi)$
is a straight line, joining the origin with $\chi(\xi, T(\xi))$.  Hence
$$\I^\alpha(\mu)~\doteq~\inf_\chi \E^\alpha(\chi)~=~\int_\Theta |\chi(\xi, T(\xi))|\, d\xi~=~\int |x|\, d\mu\,.$$

On the other hand, when $\alpha<1$, moving along a path which is traveled by few other particles
comes at a  high cost. Indeed, in this case the factor $\bigl|\chi(\xi,t)
\bigr|_\chi^{\alpha-1}$ becomes  large.   To reduce the total cost,  is thus convenient
that many particles travel along the same path.
}\end{remark}
\v
For the basic theory of ramified transport we refer 
to 
the monograph \cite{BCM}. 
For future use, we recall that optimal irrigation plans
satisfy
\v
{\bf Single Path Property:}  {\it If $\chi(\xi, \tau)=\chi(\xi',\tau')$ for some 
$\xi, \xi'\in\Theta$ and 
$0<\tau\leq \tau'$, then 
\bel{SPP}\chi(\xi, t)~=~\chi(\xi', t)\qquad\forall t\in [0, \tau].\eeq
}

\v
\subsection{The general optimization problem for branches.}
Combining the two functionals (\ref{SS2}) and (\ref{Idef}), 
one can formulate an optimization problem for the shape of branches:
\v
\begi
\item[{\bf (OPB)}] Given  a light intensity function $\eta\in \L^1(S^{d-1})$
and two constants $c>0$, $\alpha\in [0,1]$, find a positive measure $\mu$
supported on $R^d_+$ that maximizes the payoff
\bel{poff}\S^\eta(\mu)-c\,\I^\alpha(\mu).\eeq
\endi

\v
\subsection{Optimal branches in dimension $d=2$.}

\begin{figure}[ht]
	\centerline{\hbox{\includegraphics[width=8cm]{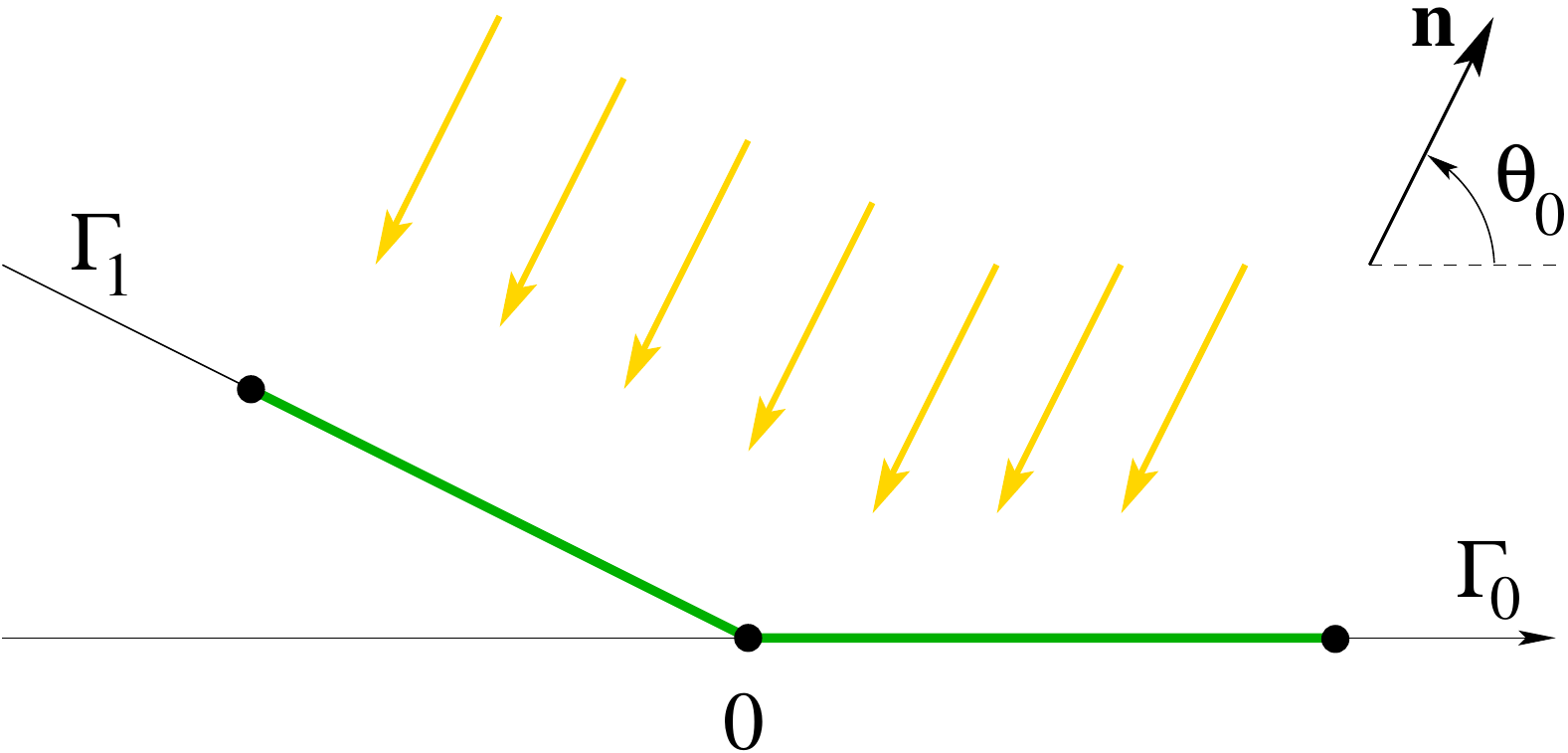}}}
	\caption{\small When the light rays impinge from a fixed direction $\bfn$,  the optimal distribution of leaves is supported on the two 
rays $\Gamma_0$ and $\Gamma_1$. } \label{f:ir100}
\end{figure}

We consider here the optimization problem for branches, in the planar case $d=2$.
 We assume  that the sunlight comes from 
a single direction $\bfn= (\cos\theta_0,
 \sin\theta_0)$,  so that the sunlight functional takes the form (\ref{SSn}).
 Moreover, as irrigation cost  we take (\ref{Idef}), for some fixed $\alpha\in \, ]0,1]$.
For a given constant $c>0$, this leads to the problem
\bel{maxb}
\hbox{maximize:}\quad \S^\bfn(\mu) -c\I^\alpha(\mu),\eeq
over all positive measures $\mu$ supported on the half space $\R^2_+\doteq
\{x=(x_1,x_2)\,;~~x_2\geq 0\}$.
To fix the ideas, we shall assume that $0<\theta_0<\pi/2$.
Our main goal is to prove that for this problem the ``solar panel" configuration shown in Fig.~\ref{f:ir100}
is optimal, namely:

\begin{theorem}\label{t:1} Assume that  $0<\theta_0\leq \pi/2$ and $1/2 \leq \alpha \leq 1$.Then the optimization problem (\ref{maxb})
has a unique solution.  The optimal measure is supported along two rays, namely
\bel{supm}
\hbox{\rm Supp}(\mu)~\subset~\Big\{ (r\cos\theta, r\sin\theta)\,;~~r\geq 0,~~
\hbox{either}~\theta=0 ~\hbox{or}~\theta =\theta_0 +{\pi\over 2}\Big\}~
\doteq~\Gamma_0\cup\Gamma_1\,.\eeq
When $0<\alpha<1/2$, the same conclusion holds provided that the angle $\theta_0$ satisfies
\bel{bigan} \cos\left( {\pi\over 2}-\theta_0\right) ~\geq~1-2^{2\alpha-1} .
\eeq
\end{theorem} 
\v

In the case $\alpha=1$ the result is straightforward.   
Indeed, for any measure $\mu$ we 
can consider its projection $\Tilde \mu$ on $\Gamma_0\cup\Gamma_1$,
obtained by shifting the mass in the direction parallel to the vector $\bfn$.   
In other words, for $x\in \R^2$ call $\phi^\bfn(x)$ the unique point in 
$\Gamma_0\cup\Gamma_1$ such that $\phi^\bfn(x)-x$ is parallel to $\bfn$.
Then let $\Tilde\mu$ be the push-forward of the measure $\mu$ w.r.t.~$\phi^\bfn$.
Since this projection satisfies $|\phi^\bfn(x)|\leq |x|$ for every $x\in \R^2_+$, 
the transportation cost decreases. On the other hand, by Remark~\ref{r:22}
the sunlight captured remains the same.  
We conclude that 
$$\S^\bfn(\Tilde \mu) - \I^1(\Tilde\mu)~\geq~\S^\bfn( \mu) - \I^1(\mu),$$
with strict inequality if $\mu$ is not supported on $\Gamma_0\cup\Gamma_1$.

In the case $0<\alpha<1$, the result is not so obvious. 
Indeed, we do not expect that the conclusion
holds if the hypothesis (\ref{bigan}) is removed. A proof of Theorem~\ref{t:1}
will be worked out in Sections 3 and 4.

Having proved that  the optimal measure $\mu$ is supported on the two rays
$\Gamma_0\cup\Gamma_1$, the density 
of $\mu$ w.r.t.~one-dimensional measure can then be determined 
using the necessary conditions derived in \cite{BGRR}.
Indeed, the density $u_1$ of $\mu$ along the ray $\Gamma_1$ provides a solution
to the scalar optimization problem
\bel{max1}
\hbox{maximize:}~~\J_1(u)~\doteq~\int_0^{+\infty} \bigl(1-e^{-u(s)}\bigr) \, ds - c \int_0^{+\infty}
\left(\int_s^{+\infty} u(r)\, dr\right)^\alpha ds\,,\eeq
among all non-negative functions $u:\R_+\mapsto\R_+$.
Here  $s$ is the arc-length variable along $\Gamma_1$.
Similarly, the density
$u_0$ of $\mu$ along the ray $\Gamma_0$ provides a solution
to the problem
\bel{max0}
\hbox{maximize:}~~\J_0(u)~\doteq~ \int_0^{+\infty} \sin\theta_0\,\bigl(1-e^{-u(s)/\sin\theta_0}\bigr) \, ds - c \int_0^{+\infty}
\left(\int_s^{+\infty} u(r)\, dr\right)^\alpha ds\,.\eeq

We write (\ref{max1}) in the form
\bel{max11}
\hbox{maximize:}~~\J_1(u)~\doteq~\int_0^{+\infty} \Big[\bigl(1-e^{-u(s)}\bigr) -
c z^\alpha \Big]ds\,,\eeq
subject to
\bel{zdot}\dot z~=~-u, \qquad z(+\infty)\,=\,0.\eeq
The necessary conditions for optimality (see for example  \cite{BP, Cesari}) 
now
yield
\bel{us}u(s)~=~\argmax_{\omega\geq 0} \Big\{- e^{-\omega} -
\,\omega q(s)\Big\}~=~-\ln q(s),\eeq
where the dual variable $q$ satisfies
\bel{qdot}\dot q~=~c\alpha z^{\alpha-1},\qquad \qquad q(0)=0.\eeq
Notice that, by (\ref{us}),  $u>0$ only if $q<1$. 
Combining (\ref{zdot}) with (\ref{qdot}) one obtains an ODE for the 
function $q\mapsto z(q)$, with  $q \in [0,1]$.  Namely
\bel{dzq} {dz(q) \over dq}~ =~ {z^{1-\alpha}  \ln{q} \over c \alpha }, \qquad
	\qquad z(1) = 0. \eeq
This equation admits the explicit solution
\bel{zq} z(q) ~= ~c^{-1/\alpha} \left[ 1 + q \ln{q} - q \right]^{1/\alpha}. \eeq
Inserting (\ref{zq}) in (\ref{qdot}), we obtain an implicit equation for 
$q(s)$:
\bel{qs} s~=~ {1 \over \alpha c^{1/\alpha} }\int_0^{q(s)}
	 \left[ 1 + t \ln{t} - t \right]^{1-\alpha \over \alpha} dt. \eeq
In turn, the density $u(s)$ of the optimal 
measure $\mu$ along $\Gamma_1$, as a function
of the arc-length $s$, is recovered from (\ref{us}).
Notice that this measure is supported only on an initial interval $[0,\ell_1]$, 	 
determined by
$$\ell_1~=~{1 \over \alpha c^{1/\alpha} }\int_0^1
	 \left[ 1 + s \ln{s} - s \right]^{1-\alpha \over \alpha} ds. $$

The density of the optimal measure along the ray $\Gamma_0$ is computed
in an entirely similar way.  In this case, the equations (\ref{us}) and  (\ref{qs}) are replaced respectively by
$$ u(s) ~=~ -(\sin \theta_0) \ln q(s),  $$
$$s~=~ {(\sin \theta_0)^{1-\alpha \over \alpha} \over \alpha c^{1/\alpha} }
\int_0^{ q(s)}
	 \left[ 1 + t \ln{t} - t \right]^{1-\alpha \over \alpha} dt.$$
Again, the condition $u(s)>0$ implies $q(s)<1$.
Along $\Gamma_0$, the optimal measure $\mu$ is supported on an initial
interval $[0, \ell_0]$, where
$$\ell_0~=~{(\sin \theta_0)^{1-\alpha \over \alpha}  \over \alpha c^{1/\alpha} }\int_0^1
	 \left[ 1 + s \ln{s} - s \right]^{1-\alpha \over \alpha} ds. $$	
	 
	 \subsection{The case $\alpha=0$.}  

In the analysis of the optimization problem {\bf (OPB)}, 
the case $\alpha =0$ stands apart.  Indeed, the general theorem on the 
existence of an optimal shape proved in \cite{BPS} does not cover this case.

When $\alpha=0$,  
a measure $\mu$ is irrigable only if it is concentrated on a set
of dimension $\leq 1$. When this happens, in any dimension $d\geq 3$ we 
have $\S^\eta(\mu)=0$
and the optimization problem is trivial.  
The only case of interest occurs in dimension $d=2$. 
In the following, $\langle\cdot, \cdot\rangle$ denotes the inner product in $\R^2$.
\v
\begin{theorem}\label{t:2}  Let $\alpha=0$, $d=2$.
Let $\eta\in \L^1(S^1)$ and 
define
\bel{a9}K~\doteq~\max_{|\bfw| = 1} ~ \int_{\bfn\in S^1} 
\Big|\langle \bfw, \bfn\rangle\Big|\, \eta(\bfn)\, d\bfn.\eeq
\begi
\item[(i)] If $K>c$, then  the optimization problem 
{\bf (OPB)}
 has no solution, because the supremum of all
possible payoffs is $+\infty$.
\item[(ii)] If $K\leq c$, then the maximum payoff  is zero,  which is trivially
achieved by the zero measure.
\endi
\end{theorem}

A proof will be given in Section~5.

\section{Properties of optimal branch configurations}
\label{s:3}
\setcounter{equation}{0}

In this section we consider the optimization problem (\ref{maxb}) in dimension 
$d=2$.  As a step toward the proof of Theorem~\ref{t:1},
some properties of  optimal branch configurations will be derived. 

 By the result in \cite{BPS} we know that an optimal measure 
$\mu$ exists and has bounded support, contained in 
$\R^2_+~\doteq~\{(x_1,x_2)\,;~~x_2\geq 0\}$.
Call $M=\mu(\R^2_+)$ the total mass of $\mu$ and let $\chi:[0,M]\times\R_+\mapsto\R^2_+$ 
be an optimal
irrigation plan for $\mu$.

\begin{figure}[ht]
	\centerline{\hbox{\includegraphics[width=8cm]{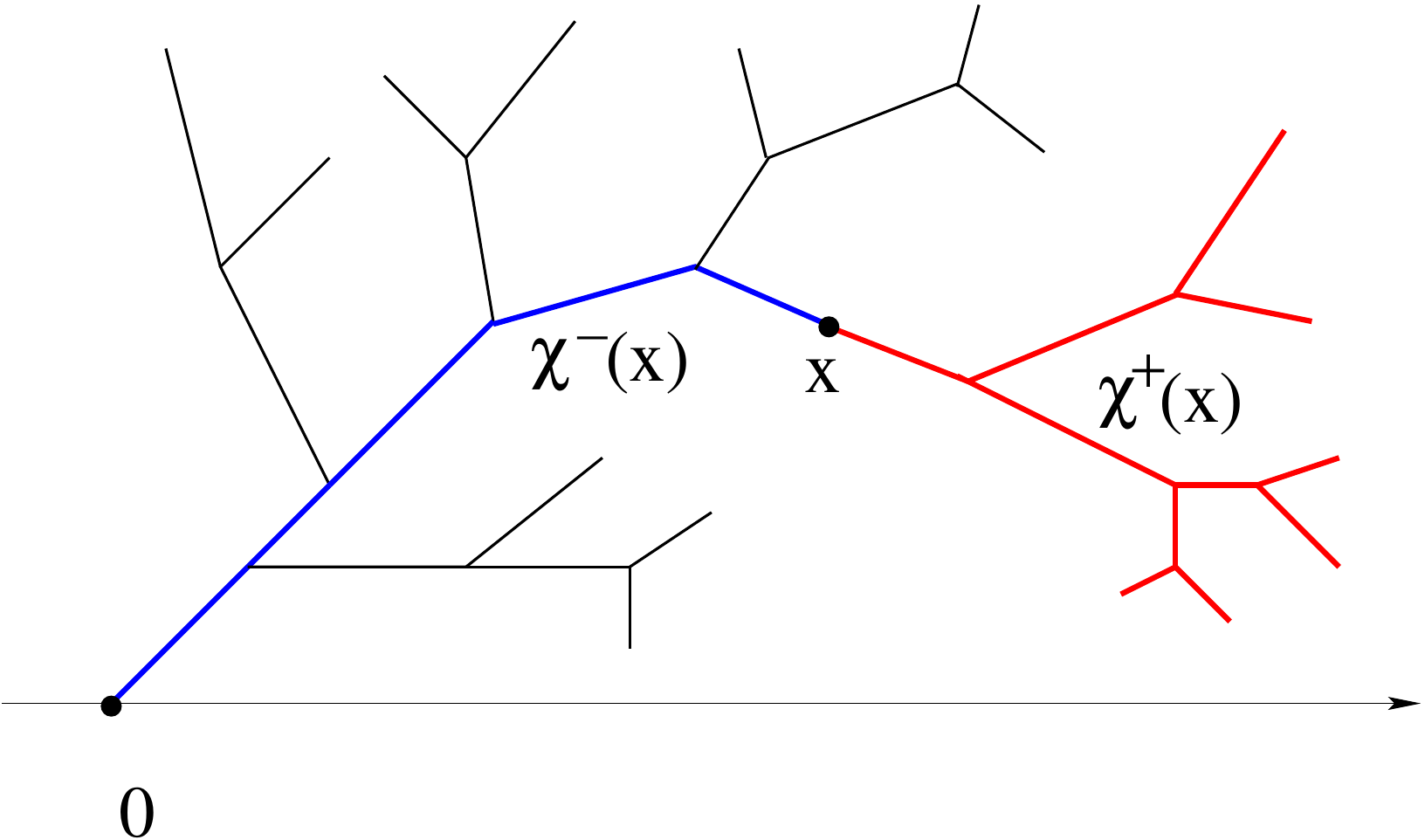}}}
	\caption{\small According to the definition~(\ref{cpm}), 
	the set $\chi^-(x)$ is a curve joining the origin to the point $x$.  The set $\chi^+(x)$
	is a subtree, containing all paths that start from $x$.} 
\label{f:ir115}
\end{figure}

Next, consider the set of all branches, namely
\bel{brs}\B~\doteq~\{ x\in\R^2_+\,;~~|x|_\chi>0\}.\eeq
By the single path property, we can introduce a partial ordering among points in $\B$.  Namely, 
for any $x,y\in \B$ we say that
$x\preceq y$ if for any $\xi\in [0,M]$ we have the implication
\bel{prec}\chi(t,\xi)\,=\,y\qquad\implies\qquad \chi(t',\xi)\,=\,x\qquad\hbox{for some}~~t'\in [0,t].\eeq
This means that all particles that reach the point $y$ pass through $x$ before 
getting to $y$.

For a given $x\in \B$ the subsets of points $y\in \B$ that precede or follow $x$ are defined as
\bel{cpm}\chi^-(x)~\doteq~\{ y\in \B\,;~~y\preceq x\},\qquad\qquad \chi^+(x)~\doteq~\{y \in \B\,;~~x\preceq y\},\eeq
respectively (see Fig.~\ref{f:ir115}).

\begin{figure}[ht]
	\centerline{\hbox{\includegraphics[width=12cm]{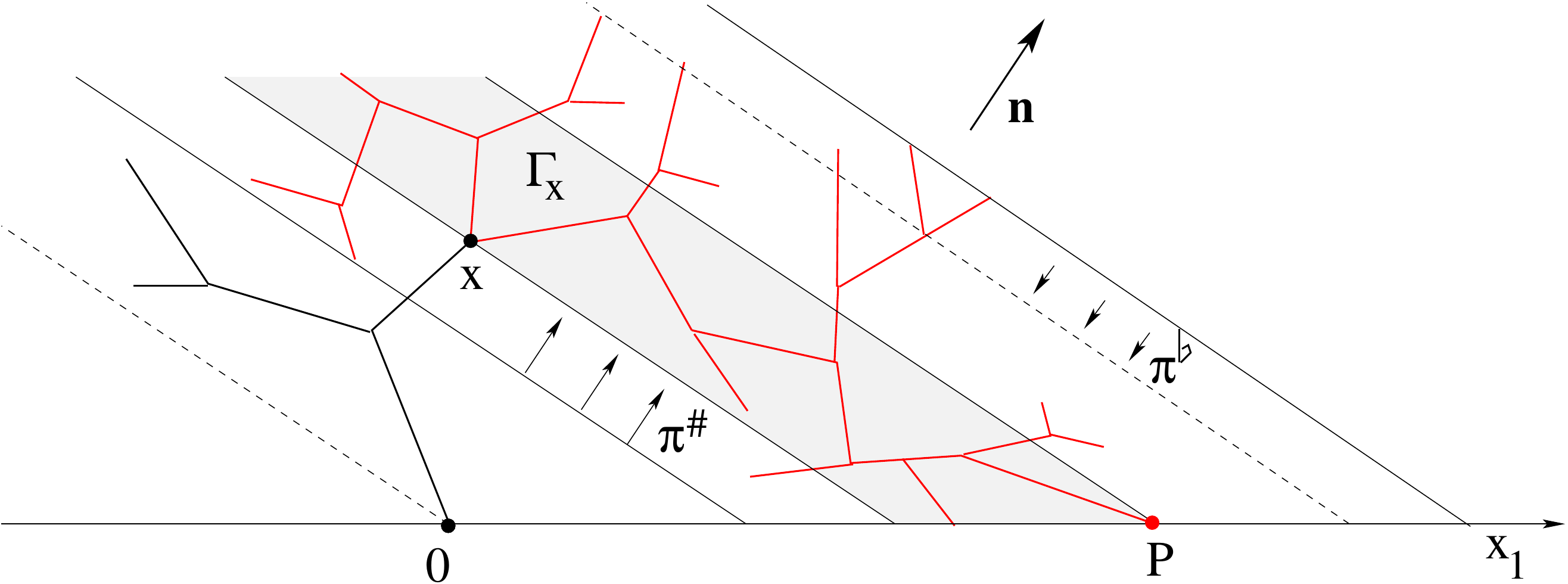}}}
	\caption{\small If the set $\chi^+(x)$ is not contained in the slab $\Gamma_x$ (the shaded region), 
	 by taking the
	perpendicular projections $\pi^\sharp$ and $\pi^\flat$ we obtain another irrigation plan with strictly lower cost,
	which irrigates a new measure  $\Tilde\mu$ gathering exactly the same amount of sunlight.
	Notice that here $P$ is the point in the closed  set $\ov{\chi^+(x)}\cap\R\bfe_1$ which has the
	largest inner product with $\bfn$. } 
\label{f:ir116}
\end{figure}

\v
We begin by deriving some properties of the sets $\chi^+(x)$.
Introducing the unit vectors
$\bfe_1 = (1,0)$, $\bfe_2=(0,1)$, we denote by $\R\bfe_1$  the set of points on the $x_1$-axis.
As before,  $\bfn= (\cos\theta_0, \sin\theta_0)$ denotes the unit vector in the direction of the sunlight.
Throughout the following, the closure of a set $A$ is denoted by $\ov A$, while
$\langle \cdot, \cdot\rangle$ denotes an inner product.
\begin{lemma}\label{l:2}
Let the measure $\mu$ provide an optimal solution to the problem (\ref{maxb}), and let 
$\chi$ be an optimal irrigation plan for $\mu$.  Then, for every $x\in \B$, one has
\bel{slab} \chi^+(x)~\subset~
\Gamma_x~\doteq~\Big\{ y\in \R^2_+\,;~~ \langle \bfn, y\rangle \,\in \,[a_x, b_x]\Big\},\eeq
where 
$a_x~\doteq~ \langle \bfn, x\rangle$, while $b_x$ is defined as follows.

\begi
\item If $\ov {\chi^+(x)}\cap \R\bfe_1=\emptyset$, then $b_x =  a_x=\langle \bfn, x\rangle$.
\item  If $\ov {\chi^+(x)}\cap \R\bfe_1\not=\emptyset$, then
$$b_x~=~\max~\{ a_x, b'_x\},\qquad\qquad b'_x~\doteq~\sup\Big\{ \langle \bfn, z\rangle\,;~~
z\in \ov{\chi^+(x)}\cap \R\bfe_1\Big\}.$$
\endi
\end{lemma}

{\bf Proof.} The right hand side of (\ref{slab}) is illustrated in Fig.~\ref{f:ir116}.
To prove the lemma, consider the set of all particles that pass through $x$, namely
$$\Theta_x~\doteq~\bigl\{ \xi\in [0,M]\,;~~\chi(\tau,\xi)=x~~\hbox{for some }~\tau\geq 0\bigr\}.$$
\v
{\bf 1.} 
We first show that, by the optimality of the solution,
\bel{bsu} \langle \bfn\,,\,\chi(\xi,t)\rangle~\geq~a_x\qquad\qquad\forall ~\xi\in \Theta_x\,,~t\geq\tau.\eeq
Indeed, consider the perpendicular
projection on the half plane
$$\pi^\sharp:\R^2~\mapsto~ S^\sharp
~\doteq~\{y\in\R^2\,;~\langle \bfn, y\rangle~\geq~a_x\}.$$
Define the projected irrigation plan
$$\chi^\sharp(t,\xi)~\doteq~\left\{ \bega{rl} \pi^\sharp\circ \chi(t,\xi)\qquad &\hbox{if} ~~\xi\in \Theta_x\,,~~t\geq \tau,
\\[3mm]  \chi(t,\xi)\qquad &\hbox{otherwise.}\enda\right.
$$
Then the new measure $\mu^\sharp$ irrigated by $\chi^\sharp$ is still supported on $\R^2_+$ and 
has exactly the same projection on 
$E^\perp_\bfn$ as $\mu$. Hence it gathers the same amount of sunlight.  However,
if the two irrigation plans do not coincide a.e., then the cost of $\chi^\sharp$ is strictly smaller than
the cost of $\chi$,
contradicting the optimality assumption.
\v
{\bf 2.}
Next, we show that 
\bel{bsw} \langle \bfn\,,\,\chi(\xi,t)\rangle~\leq~b_x\qquad\qquad\forall ~\xi\in \Theta_x\,~t\geq\tau.\eeq
Indeed, call
$$b''~\doteq~\sup~\Big\{ \langle \bfn, z\rangle\,;~~
z\in \chi^+(x)\Big\}.$$
If $b''\leq b_x$, we are done.   In the opposite case, by a continuity and compactness argument
we can find $\delta>0$ such that the following holds.
Introducing the perpendicular
projection on the half plane
$$\pi^\flat:\R^2~\mapsto~ S^\flat
~\doteq~\{y\in\R^2\,;~\langle \bfn, y\rangle~\leq~b''-\delta\},$$
one has
\bel{ag}\bigl\{ \pi^\flat(y)\,;~~y\in \chi^+(x)\bigr\}
~\subseteq~\R^2_+\,.\eeq
Similarly as before, define the projected irrigation plan
$$\chi^\flat(t,\xi)~\doteq~\left\{ \bega{rl} \pi^\flat\circ \chi(t,\xi)\qquad &\hbox{if} ~~\xi\in \Theta_x\,,~~t\geq \tau,
\\[3mm]  \chi(t,\xi)\qquad &\hbox{otherwise.}\enda\right.
$$
Then the new measure $\mu^\flat$ irrigated by $\chi^\flat$  is supported on $\R^2_+\cap S^\flat$ and has exactly the same projection on 
$E^\perp_\bfn$ as $\mu$. Hence it gathers the same amount of sunlight.  However,
if the two irrigation plans do not coincide a.e., then the cost of $\chi^\flat$ is strictly smaller than
the cost of $\chi$,
contradicting the optimality assumption.
This completes the proof of the Lemma.
\endproof

\begin{figure}[ht]
	\centerline{\hbox{\includegraphics[width=10cm]{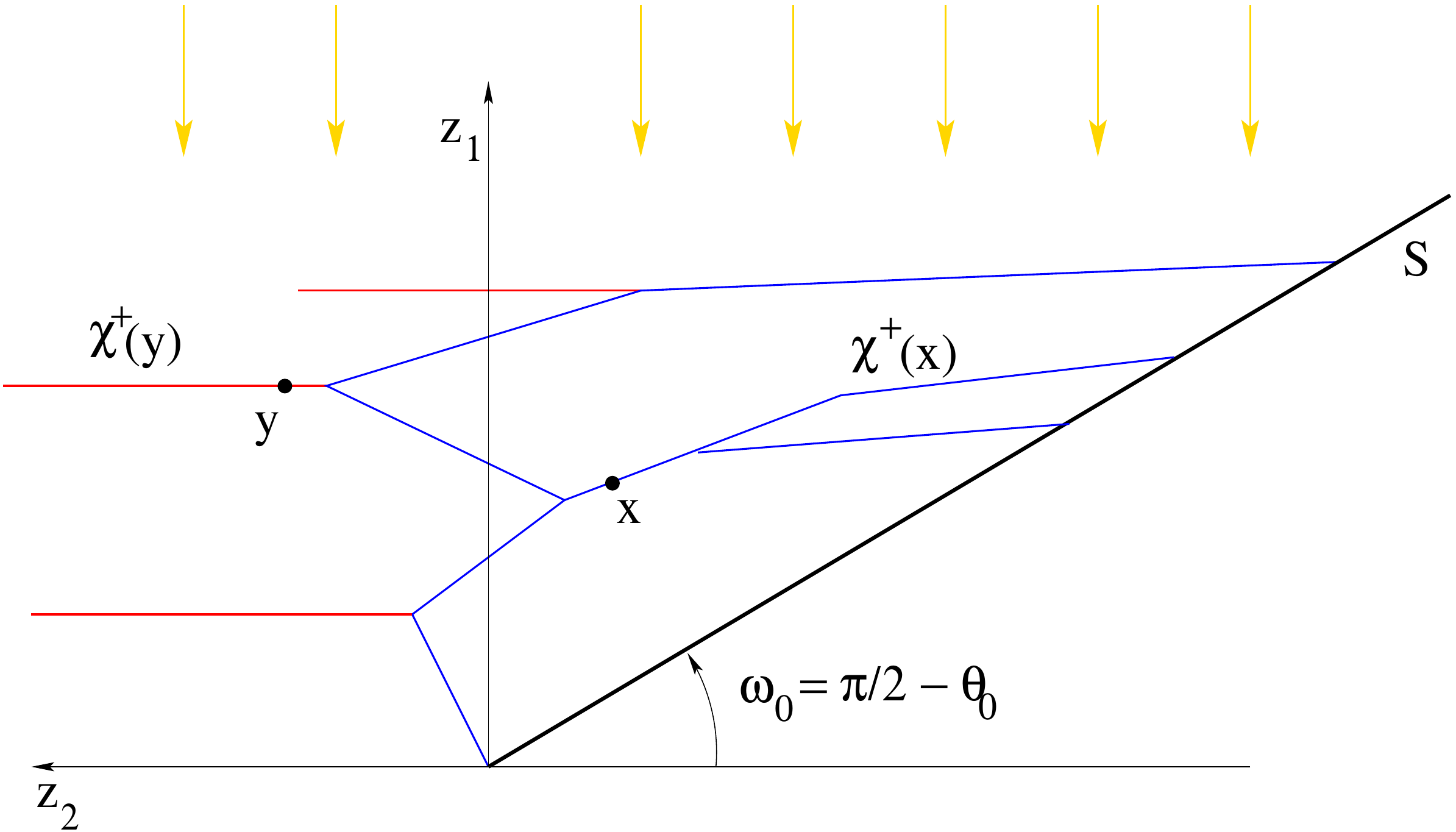}}}
	\caption{\small  After a rotation of coordinates, the sunlight comes from the vertical direction.
	Here the blue lines correspond to the set $\B^*$ in \eqref{B*}.} 
\label{f:ir123}
\end{figure}

\v
 Based on the previous lemma, we now consider the set
\bel{B*}\B^*~\doteq~\{ x\in \B\,;~~\ov{\chi^+(x)}\cap\R\bfe_1\not=\emptyset\}.\eeq
It will be convenient to rotate coordinates by an angle of $\pi/2 - \theta_0$, and choose
new coordinates $(z_1, z_2)$ oriented as in Fig.~\ref{f:ir123}.
In these new coordinates, the direction of sunlight becomes vertical, while
the positive $x_1$-axis corresponds to the line
\bel{bfS} \bfS~\doteq~\bigl\{(z_1, z_2)\,;~~z_1\geq 0\,,\quad z_2 = 
-\lambda z_1\bigr\},
\qquad\hbox{where}\quad   \lambda =  \tan\theta_0\,.\eeq

Calling $\bigl(z_1(\xi,t), z_2(\xi,t)\bigr)$ the corresponding coordinates of the point $\chi(\xi,t)$,
from Lemma~\ref{l:2} we immediately obtain

\begin{corollary}\label{c:2}  Let $\chi$ be an optimal irrigation plan for a solution to (\ref{maxb}). Then
\begi
\item[(i)] For every $\xi\in [0,M]$, the map $t\mapsto z_1(\xi,t)$ is non-decreasing.
\item[(ii)] If $\bar z=(\bar z_1,\bar z_2)\notin \B^*$, then $\chi^+(\bar z)$
is contained in a horizontal line.  Namely, 
\bel{hor}\chi^+(\bar z)\subset \{ (\bar z_1, s)\,;~~s\in\R\}.\eeq
\endi
\end{corollary}

To make further progress, we define
 $$z_1^{\max}~\doteq~\sup\,\bigl\{ z_1\,;~(z_1,\,z_2)\in \B^*\bigr\}.$$ 
Moreover, on the interval $[0, z_1^{max}[\,$ we consider the function
\bel{phid}
\vp(z_1)~\doteq~\sup\,\bigl\{ s\,;~~(z_1,s)\,\in\,\B^*\bigr\}.\eeq

\begin{figure}[ht]
	\centerline{\hbox{\includegraphics[width=12cm]{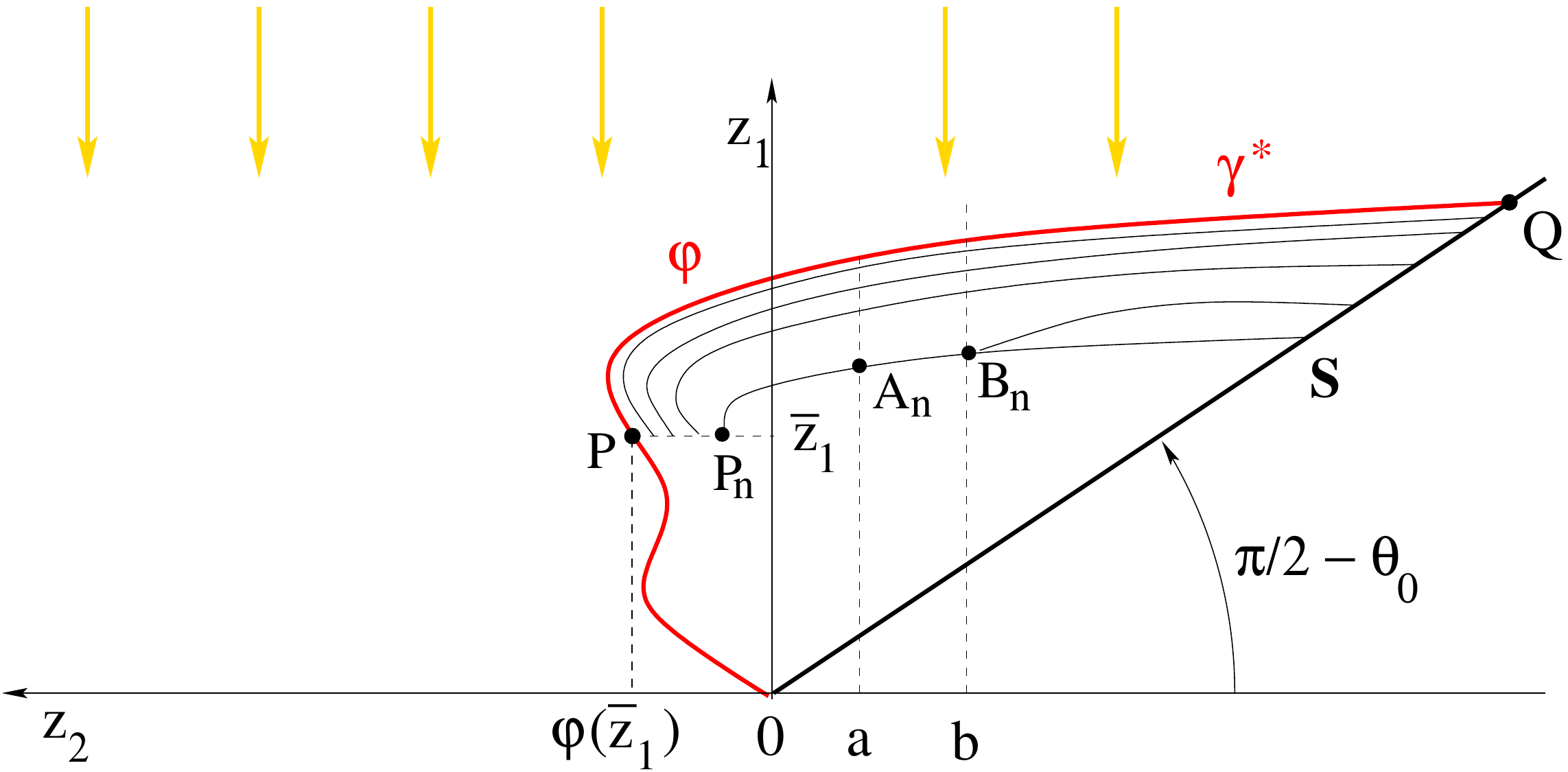}}}
	\caption{\small  The construction used in the proof of Lemma~\ref{l:33}. } 
\label{f:ir125}
\end{figure}

\begin{lemma}\label{l:33}
For every $z_1\in [0, z_1^{max}[\,$, the supremum $\vp(z_1)$ is attained as a maximum.  
\end{lemma}

{\bf Proof.} 
{\bf 1.} Assume that, on the contrary, for some $\bar z_1$  the supremum is not a maximum.
In this case,  as shown in Fig.~\ref{f:ir125}, there exist a sequence of points 
$P_n\to P$ with $P_n = (\bar z_1, s_n)$, $P=(\bar z_1, \bar z_2)$, $s_n\uparrow \bar z_2$.
Here $P_n\in \B^*$ for every $n\geq 1$ but $P\notin \B^*$.
\v
{\bf 2.} Choose  two values $a,b$ such that 
$$-\lambda \bar z_1~ <~b~<~a~<~\vp(\bar z_1).$$ 
By construction, for every $n\geq 1$ the set $\ov {\chi^+(P_n)}$ intersects 
 $\bfS$.
Therefore we can find points
$$P_n~\prec~A_n~\prec~B_n$$
all in $\B^*$, with 
$$A_n~=~(t_n, a),\qquad B_n~=~(t'_n, b),
\qquad\qquad \bar z_1\leq t_n \leq t'_n \leq z_1^{max}\,.$$
\v
{\bf 3.}
Since the branches $\chi^+(A_n)$ are all disjoint, we have
$$\sum_{n\geq 1} |A_n|_\chi~\leq~M~\doteq~\mu(\R^2_+).$$
We can thus find $N$ large enough so that 
\bel{eN}\ve_N~\doteq~|A_N|_\chi~< ~(a-b)^{1\over 1-\alpha}.\eeq

Consider the modified transport plan $\Tilde\chi$,  obtained 
from $\chi$ by removing all particles that go through the point $B_N$.
More precisely, 
$\Tilde\chi$ is the restriction of $\chi$ to the domain
$$\Tilde\Theta~\doteq~\Theta\setminus \{ \xi\,;~~\chi(\xi,\tau)=B_N\quad
\hbox{for some}~\tau\geq 0\}.$$
Let $\Tilde\mu$ be the measure irrigated by $\Tilde\chi$.

Since $\Tilde\mu\leq\mu$, the total amount of sunlight 
gathered by the measure $\Tilde\mu$ satisfies
\bel{e1}\S(\mu)-\S(\Tilde \mu)~\leq~(\mu-\Tilde\mu)(\R^2) .\eeq

We now estimate the reduction in the transportation cost, achieved by replacing 
$\mu$ with $\Tilde \mu$.   Since all water particles reaching $B_N$ must pass through
$A_N$, they must cover a distance $\geq |B_N-A_N|\geq a-b$ traveling 
along a path whose maximum flux is $\leq \ve_N$.
The difference in the  transportation  costs can thus be estimated by
\bel{e2}\I^\alpha( \mu)-\I^\alpha(\Tilde \mu)~\geq~
(a-b)\cdot \ve_N^{\alpha-1}\cdot (\mu-\Tilde\mu)(\R^2).\eeq
If (\ref{eN}) holds, combining (\ref{e1})-(\ref{e2}) we obtain
$$\S(\mu) -\I^\alpha( \mu)~<~\S(\Tilde\mu) -\I^\alpha( \Tilde\mu).$$
Hence the measure $\mu$ is not optimal.  This contradiction proves the lemma.
\endproof

\begin{figure}[ht]
	\centerline{\hbox{\includegraphics[width=12cm]{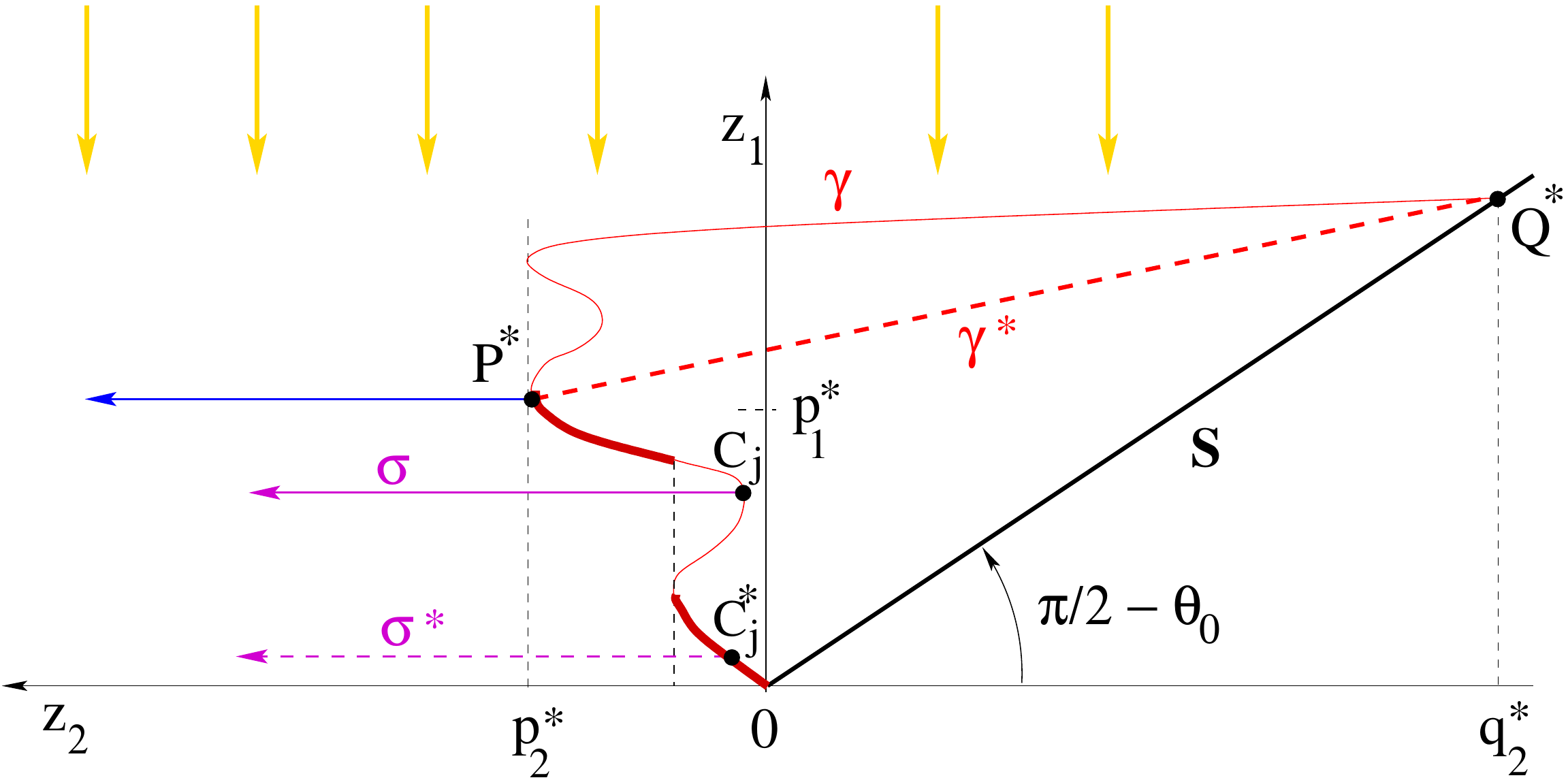}}}
	\caption{\small The thick portions of the curve $\gamma$ are the only
	points where a left bifurcation can occur. 
	If a horizontal branch $\sigma$ bifurcates from $C_j$, all the mass on this 
	branch can be shifted downward to another branch $\sigma^*$ 
	bifurcating from $C_j^*$.
Furthermore, if some portion of the path $\gamma$ between $P^*$ and $Q$
lies above the segment $\gamma^*$ joining these two points, we can take a projection of $\gamma$ on $\gamma^*$.
In both cases,   the transportation cost is strictly reduced.
} 
\label{f:ir126}
\end{figure}

\v
By the previous result, the graph of $\vp$ is contained in one single maximal trajectory of the transport plan $\chi$.   As in Figure~\ref{f:ir126}, we 
let $s\mapsto \gamma(s)$ be the arc-length parameterization of this curve,  
which provides the 
left boundary of the set $\B^*$.  
 
Along the curve $\gamma$, we now 
consider the set of points $C_j = (z_{1,j}, z_{2,j})$ where some horizontal branch
bifurcates on the left.   A property of such points is given below.

\begin{lemma}\label{l:34}
In the above setting, for every $j$, one has
\bel{phij}
\vp(s)~<~z_{2,j}\qquad\forall s<z_{1,j}\,.\eeq
\end{lemma}

{\bf Proof.} 
If (\ref{phij}) fails, there exists another point $C_j^*= (z_{1,j}^*,z_{2,j})$ 
along the curve $\gamma$, with $z_{1,j}^*<z_{1,j}$. 
We can now replace the measure $\mu$ by another measure $\Tilde\mu$ obtained 
as follows.   All the mass lying on the horizontal half-line 
$\{(z_{1,j},s)\,;~~s\geq z_{2,j}\}$ is shifted downward on the half-line 
$\{(z_{1,j}^*,s)\,;~~s\geq z_{2,j}\}$.  
Since the functional $\S^\bfn$ is invariant under 
 vertical shifts, we have
 $\S^\bfn(\Tilde\mu) = \S^{\bfn}(\mu)$. However, the transportation cost
 is strictly smaller:    $\I^\alpha(\Tilde\mu) < \I^\alpha(\mu)$.  This contradicts the optimality of $\mu$. \endproof
\v

Next, as shown in Fig.~\ref{f:ir126}, we consider a point $P^*= (p_1^*, p_2^*)\in \gamma$ where the component
$z_2$ achieves its maximum, namely
\bel{z22}
p_2^*~=~\max\{z_2\,;~~(z_1, z_2)\in \gamma\}~\geq~0.\eeq
Notice that such a maximum exists because $\gamma$ is a continuous curve,
starting at the origin. If this maximum is attained at more than one point, 
we choose the one with smallest $z_1$-coordinate, so that
\bel{z11} p_1^* ~=~\min\{z_1\,;~~(z_1, p_2^*)\in \gamma\}.\eeq

Moreover, call 
$$q_2^*~\doteq~\inf\{ z_2\,;~~(z_1,z_2)\in \hbox{Supp}(\mu)\},$$
and let $Q^*= (q_1^*, q_2^*) \in \bfS$ be the point on the ray $\bfS$ whose
second coordinate is $p_2^*$.  We observe that, by the optimality of the solution,
all paths of the irrigation plan $\chi$ must lie within the convex set
$$\Sigma^*~\doteq~\{(z_1,z_2)\,;~~z_1\in [0, q_1^*],\quad z_2\geq q_2^*\}.$$
Otherwise, calling $\pi^*:\R^2\mapsto \Sigma^*$ the perpendicular projection
on the convex set $\Sigma^*$, the composed plan
$$\chi^*(\xi,t)~\doteq~\pi^*\bigl(\chi(\xi,t)\bigr)$$
would satisfy
$$\S^\bfn(\chi^*)~=~\S^\bfn(\chi),\qquad \E^\alpha(\chi^*)~<~
\E^\alpha(\chi),$$
contradicting the optimality assumption.

By a  projection argument we now show  that, in an optimal solution, all the particle paths remain below the 
segment $\gamma^*$ with endpoints $P^*$ and $Q^*$.

\begin{lemma}\label{l:35}
In the above setting, let 
$$\gamma^*~=~\bigl\{(z_1,z_2)\,;~z_1 = a+bz_2\,,\qquad z_2\in [q_2^*, p_2^*]
\bigr\}$$
be the segment with endpoints $P^*, Q^*$.
If 
\bel{oip}(\xi,t)\mapsto \chi(\xi,t)~=~(z_1(\xi,t), z_2(\xi,t)\bigr)\eeq
is an optimal irrigation plan for the problem (\ref{maxb}), then 
we have the implication
\bel{below}
z_2(\xi,t)\,\in\, [q_2^*, p_2^*]\qquad\implies\qquad z_1(\xi,t)~\leq ~a+b\,z_2(\xi,t).
\eeq
\end{lemma}
\v
{\bf Proof.} {\bf 1.} 
It suffices to show that the maximal curve $\gamma$ 
lies below $\gamma^*$.    If this is not the case,
consider the set of particles which go through the point $P^*$ and 
then move to the right of $P^*$, namely
\bel{OS}
\Omega^*~=~\Big\{\xi\in [0,M]\,;~~\chi(\xi,t^*)= P^* ~~\hbox{for some $t^*\geq 0$}, 
~~~z_2(\xi, t)< p_2^*~~\hbox{for $t>t^*$}\Big\}.\eeq
\v
{\bf 2.} Consider the convex region below $\gamma^*$, defined by
$$\Sigma~\doteq~\Big\{ (z_1,z_2)\,;~~0\le z_1\leq a+ bz_2\,,\qquad z_2\in 
[q_2^*, p_2^*]\Big\}.$$
Let $\pi:\R^2\mapsto\Sigma$ be the perpendicular projection.
Then the irrigation plan
\bel{chid}\chi^\dagger(\xi,t)~\doteq~\left\{\bega{rl} \pi \Big(\chi(\xi,t)\Big)\quad &\hbox{if}~~
\xi\in \Omega^*, ~t>t^*,\\[3mm]
\chi(\xi,t)\quad &\hbox{otherwise,}\enda\right.\eeq
has total cost strictly smaller than $\chi$.
Indeed, for all $x, \xi,t$ we have
\bel{equal}\bigl|\pi(x)\bigr|_{\chi^\dagger}~\geq~|x|_\chi\,,\qquad
\bigl|\dot \chi^\dagger(\xi,t)\bigr|~\leq~\bigl|\dot \chi(\xi,t)\bigr|.\eeq
Notice that, in (\ref{equal}), equality can  hold for a.e.~$\xi,t$ 
only in the case where $\chi=\chi^\dagger$.

\v
{\bf 3.} We now observe that the perpendicular projection on $\Sigma$
can decrease the $z_2$-component.   As a consequence, the measures $\mu$ and $\mu^\dagger$  irrigated by  $\chi$ and $\chi^\dagger$ may have a different projections on the $z_2$ axis.
If this happens, we may have $\S^\bfn(\mu)\not=\S^\bfn(\mu^\dagger)$.

To address this issue, we observe that all particles $\xi\in \Omega^*$
satisfy
$\chi^\dagger(\xi, t^*)= \chi(\xi,t^*)=P^*$. In terms of the $z_1,z_2$ 
coordinates, this implies
\bel{zp2}
z_2^\dagger (\xi,t^*) ~=~ z_2(\xi, t^*) ~= ~p_2^*,
\qquad z_2^\dagger(\xi, T(\xi))~ \leq ~z_2(\xi,T(\xi)) ~< ~p_2^*\,.\eeq
By continuity, for each $\xi\in\Omega^*$  we can find a stopping
time $\tau(\xi)\in [t^*, T(\xi)]$ such that
$$z_2^\dagger(\xi, \tau (\xi))~=~z_2(\xi,T(\xi)).$$
Call $\Tilde \chi$ the truncated irrigation plan, such that
\bel{cds}\Tilde \chi(\xi,t)~\doteq~
\left\{\bega{cl} \chi^\dagger (\xi,t)\quad &\hbox{if}~~
\xi\in \Omega^*, ~t\leq \tau(\xi),\\[3mm]
\chi(\xi,\tau(\xi))\quad \quad &\hbox{if}~~
\xi\in \Omega^*, ~t\geq \tau(\xi),\\[3mm]
\chi(\xi,t)\quad &\hbox{if}~~\xi\notin\Omega^*.\enda\right.\eeq
By construction, the measures $\mu$ and $\Tilde \mu$  irrigated by  $\chi$ and $\Tilde \chi$ have exactly the same projections on the $z_2$ axis.
Hence $\S^\bfn(\Tilde\mu)~=~\S^\bfn(\mu)$. On the other hand, 
the corresponding costs satisfy
$$\E^\alpha(\Tilde\chi)~\leq~\E^\alpha(\chi^\dagger)~<~\E^\alpha(\chi).$$
This contradicts optimality, thus proving the lemma.\endproof

\section{Proof of Theorem~\ref{t:1}}
\label{s:4}
\setcounter{equation}{0}

In this section we give a proof of Theorem~\ref{t:1}.
As shown in Fig.~\ref{f:ir126},
let $P^*=(p_1^*, p_2^*)$ be the point defined at (\ref{z22})
We consider two cases:
\begi
\item[(i)] $P^*=0\in \R^2$,
\item[(ii)] $P^*\not= 0$.
\endi

Assume that  case (i) occurs. 
Then, by Lemma~\ref{l:34},
the only branch that can bifurcate to the left of $\gamma$ 
must lie on the $z_2$-axis.    Moreover, by Lemma~\ref{l:35},
the path $\gamma$  cannot lie above the segment with endpoints $P^*$, $Q^*$. Therefore,  the restriction of the measure $\mu$ to the half space $\{z_2\leq 0\}$
is supported on the  line $\bfS$.
Combining these two facts we achieve the conclusion of the theorem.

The remainder of the proof will be devoted to showing that the case (ii)
cannot occur, because it  would contradict the optimality of the 
solution.
\v

\begin{figure}[ht]
	\centerline{\hbox{\includegraphics[width=12cm]{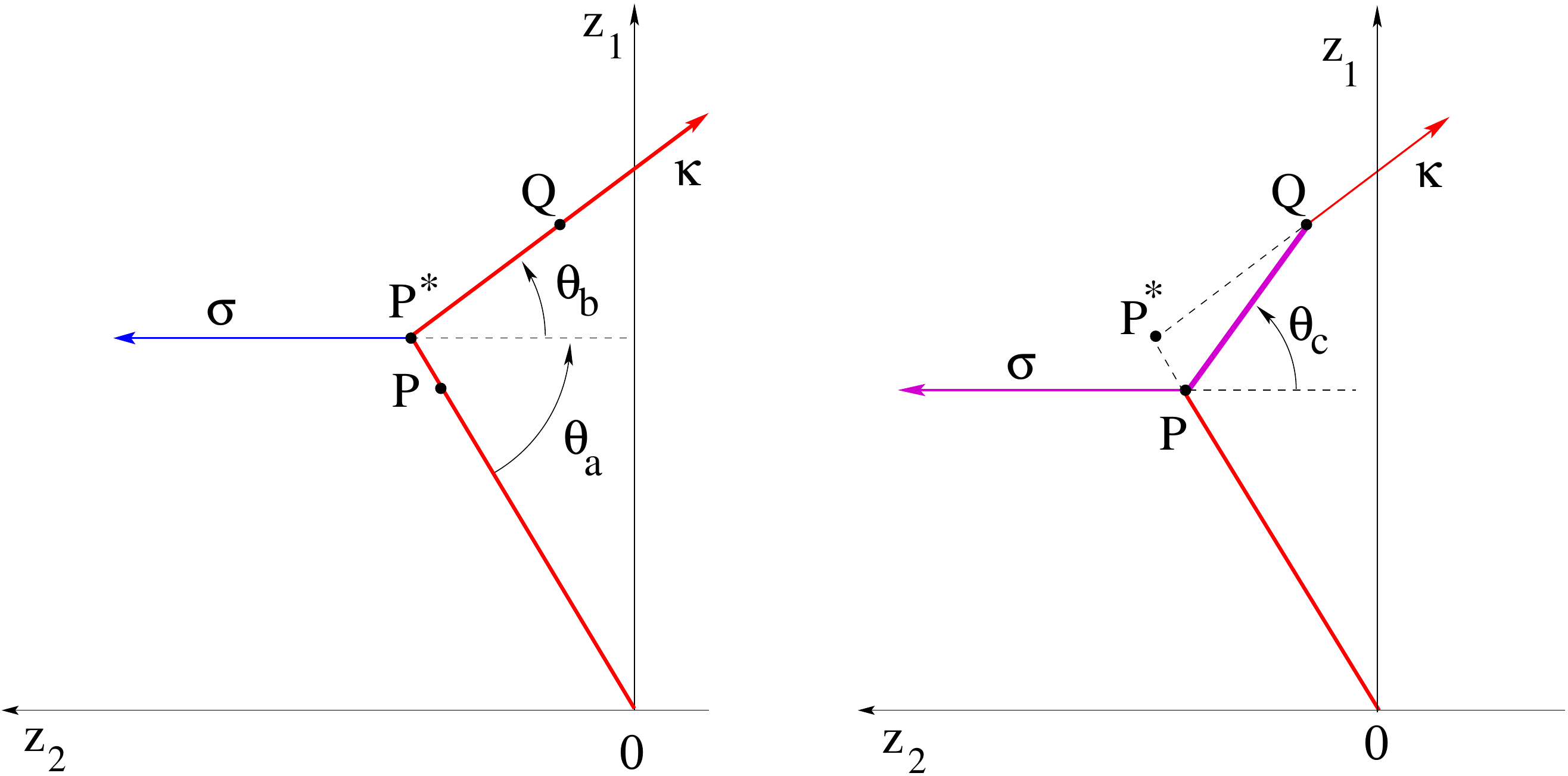}}}
	\caption{\small The basic case:  in  a neighborhood of $P^*$ the 
	trajectories are straight lines.   To show that the configuration on the left is not optimal, we replace the  portion of the trajectory 
	between $P$ and $Q$ with a single segment. } 
	\label{f:ir129}
\end{figure}

To illustrate the heart of the matter, we first consider the  elementary configuration
shown in Fig.~\ref{f:ir129}, left, where all trajectories are straight  lines. 
We call $\kappa$ the flux along the segment $P^*Q$ 
and $\sigma$ the flux along the 
horizontal line bifurcating to the left of $P^*$.
As in  Fig.~\ref{f:ir129}, right, we then 
replace the segments $P P^*$ and $P^* Q$ by a single segment with endpoints $P,Q$.
To fix the ideas, the lengths of these two segments will be denoted by
\bel{lab}\ell_a~=~|P-P^*|,\qquad \qquad \ell_b~=~|Q-P^*|.\eeq
The angles between these segments and a horizontal line will be denoted by
$\theta_a,\theta_b$, respectively.   Our main assumption is
\bel{abas} 0\,\leq \,\theta_a\,\leq\, {\pi\over 2},\qquad\quad
  0\,\leq \,\theta_b\,<\, {\pi\over 2} - \theta_0\,.
\eeq

Having performed this modification, the previous transportation cost along $PP^*$ and $P^*Q$
$$ (\kappa + \sigma)^\alpha \ell_a + \kappa^\alpha \ell_b $$
is replaced by 
\bel{new} \kappa^\alpha \sqrt{ \ell_a^2 + \ell_b^2 
- 2 \ell_a \ell_b \cos(\theta_a + \theta_b) }+ \sigma^\alpha \ell_a\cos\theta_a \,. \eeq
Notice that the last term in (\ref{new}) accounts for the fact that an amount $\sigma$ of particles
need to cover a longer horizontal distance, reaching  $P$ instead of $P^*$.

The difference in the cost is thus expressed by the function
$$ f(\ell_a,\ell_b) ~= ~(\kappa + \sigma)^\alpha \ell_a - \sigma^\alpha \ell_a\cos\theta_a 
+ \kappa^\alpha \left[ \ell_b - \sqrt{ \ell_a^2 + \ell_b^2 - 2 \ell_a \ell_b \cos(\theta_a + \theta_b) } \right]. $$
Notice that  this function is positively homogeneous of degree 1 w.r.t.~the variables 
$\ell_a,\ell_b$.
We observe that, by choosing the angle $\theta_c$ between the 
segment $PQ$ and a horizontal line to be just slightly larger than
$\theta_b$, we can render the ratio
$\ell_a/\ell_b$ as small as we like.   Taking advantage of this fact, we set
$$\ell_a=\ve\ell,\qquad \ell_b~=~\ell$$ 
for some $\ve>0$ small.
By the homogeneity of $f$ it follows
$$ f(\ve \ell, \ell)~ = ~\ell \left[ \ve (\kappa + \sigma)^\alpha - \ve \sigma^\alpha \cos\theta_a 
+ \kappa^\alpha \Big( 1 - \sqrt{ 1 + \ve^2 - 2 \ve \cos(\theta_a + \theta_b) } \Big) \right].$$
This yields
\bel{pe2} \bega{rl} \ds{d\over d\ell} f(\ve \ell, \ell) &=~ \ve (\kappa + \sigma)^\alpha - \ve \sigma^\alpha \cos\theta_a 
+ \kappa^\alpha \Big( 1 - \sqrt{ 1 + \ve^2 - 2 \ve \cos(\theta_a + \theta_b) } \Big) \\[3mm]
&= ~\ds \ve \Big[ (\kappa + \sigma)^\alpha - \sigma^\alpha \cos\theta_a + \kappa^\alpha \cos(\theta_a + \theta_b) + \O(1)\cdot \ve \Big] .\enda \eeq


Setting
$$\lambda ~=~{\sigma\over \kappa+\sigma}$$
we now study the function 
\bel{Fdef}F(\lambda, \theta_a,\theta_b)~\doteq~1-
\lambda^\alpha \cos\theta_a +(1- \lambda)^\alpha \cos(\theta_a +\theta_b ),\eeq
and find under which conditions on $\theta_b$
this function $F$ it remains positive for all $\lambda\in [0,1]$,  $\theta_a\in [0, \pi/2]$.

\begin{lemma} \label{l:41}

\begi
\item[(i)]
For $\alpha\geq 1/2$ and any $\theta_a,\theta_b\in [0, \pi/2]$, we always have 
$F(\lambda, \theta_a,\theta_b)\geq 0$.
\item[(ii)] When $0<\alpha<1/2$ we have
$F(\lambda, \theta_a,\theta_b)\geq 0$ for every $\theta_a,\theta_b\in [0,\pi/2]$
provided that $\theta_b$ satisfies the additional bound
\bel{tbb} \cos \theta_b ~\ge ~1- 2^{2\alpha-1}.\eeq
\endi
\end{lemma}
\v
{\bf Proof.}
The function $F$ in (\ref{Fdef}) can be written in terms of an inner product:
\bel{fip} \bega{rl} F(\lambda,\theta) &=~ 1 - \cos{\theta_a}\left[\lambda^\alpha - (1-\lambda)^\alpha \cos{\theta_b} \right] - \sin{\theta_a} (1-\lambda)^\alpha \sin{\theta_b} \\[4mm]
&= ~1 - \Big\langle \left(\cos{\theta_a}, \,\sin{\theta_a}\right)\,,~\Big( \lambda^\alpha - (1-\lambda)^\alpha \cos{\theta_b}~,~ (1-\lambda)^\alpha \sin{\theta_b} \Big)\Big\rangle.
 \enda\eeq
To prove that $F\geq 0$ it thus suffices to show that the second vector on the right hand side of 
(\ref{fip}) has length less than or equal to one. Namely
$$ \lambda^{2\alpha} + (1-\lambda)^{2\alpha} 
- 2\lambda^\alpha (1-\lambda)^\alpha \cos{\theta_b} ~\le ~1. $$
This inequality holds provided that
\bel{c3} \cos{\theta_b} ~\ge ~{ \lambda^{2\alpha} + (1-\lambda)^{2\alpha} -1 \over 2 \lambda^\alpha (1-\lambda)^\alpha }\,. \eeq
In the case where $\alpha\geq 1/2$ we have
$$\lambda^{2\alpha} + (1-\lambda)^{2\alpha} ~\leq~1\qquad\hbox{for all}~\lambda\in [0,1],$$ 
hence  (\ref{c3}) holds.

To study the case where $\alpha<1/2$, 
consider the function
$$ g(\lambda)~ \doteq~ {  \lambda^{2\alpha} + (1-\lambda)^{2\alpha} -1\over 2 \lambda^\alpha (1-\lambda)^\alpha }~=~1+ { \bigl( \lambda^{\alpha} - (1-\lambda)^{\alpha}\bigr)^2 -1\over  2 \lambda^\alpha (1-\lambda)^\alpha}\, . $$
We observe  that, for $0 \le \alpha \le \tfrac12$, one has
\bel{gp}0~\leq ~g(\lambda)~\leq ~g\Big({1\over 2}\Big) ~=~1-  2^{2\alpha-1} ,\eeq
while 
$$\lim_{\lambda\to 0+} g(\lambda)~=~\lim_{\lambda\to 1} g(\lambda)~=~0.$$

From (\ref{gp}) it now follows that the condition (\ref{tbb}) guarantees that 
(\ref{c3}) holds, hence $F\geq 0$, as required.
\endproof



\begin{figure}[ht]
	\centerline{\hbox{\includegraphics[width=10cm]{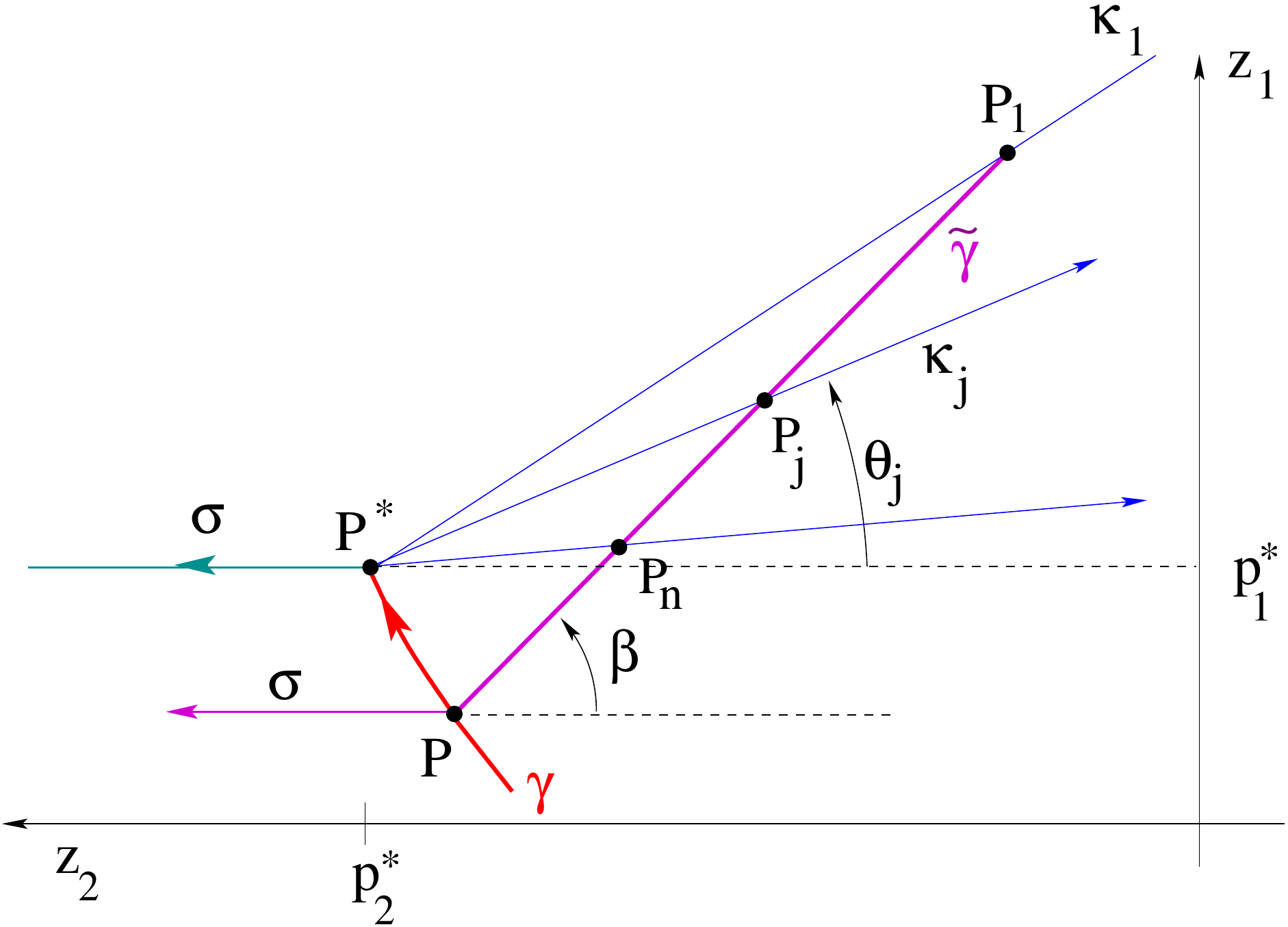}}}
	\caption{\small  A more general configuration, compared with the one in Fig.~\ref{f:ir129}.   
	} 
	\label{f:ir135}
\end{figure}

We now consider the more general configuration shown in Fig.~\ref{f:ir135}.
Water is transported along the path $\gamma$ up to the point $P^*$.
Then the flux is split into a finite number of straight paths.
One goes horizontally to the left, with flux $\sigma\geq 0$.   The other 
pipes go
 to the right, with fluxes $\kappa_1,\ldots,\kappa_n >0$, 
 at angles 
\bel{angles}
0~\leq ~\theta_n~<~\cdots~<~\theta_2~<~\theta_1~< ~{\pi\over 2}-\theta_0.\eeq

We compare this configuration with a modified irrigation plan, 
where a ``bypass" is inserted along a segment $\Tilde\gamma$ with 
endpoints $P$, $P_1$,
at an angle $\beta$ satisfying
\bel{beta} \theta_1~<~\beta~< ~{\pi\over 2}-\theta_0.\eeq
In this case, water particles travel  along $\gamma$ until they reach $P$.
Then, an amount $\sigma$ of particles bifurcates to the left. All the remaining 
particles are transported along the segment  $\Tilde\gamma$,
until they reach the points $P_n,\ldots, P_1$ along the old pipes.
The next lemma estimates the saving in the irrigation cost 
achieved by inserting the ``bypass" along the segment $PP_1$.

\begin{lemma}
\label{l:42}  As in Theorem~\ref{t:1}, assume that either
$1/2\leq\alpha \leq 1$, or else (\ref{bigan}) holds.
In the above setting, one has 
\bel{saving}
\hbox{\rm [old cost]} - \hbox{\rm [new cost]}~\geq~|P_1-P^*|\cdot \delta(\theta_1,\kappa),\eeq
where $\delta(\theta_1, \kappa)$ is a continuous function,
strictly positive
for  $0\leq \theta_1< {\pi\over 2}-\theta_0$ and $\kappa=\kappa_1+\cdots+\kappa_n>0$.
\end{lemma} 

{\bf Proof.}
 {\bf 1.}  As in the previous lemmas, we call $\theta_a$ the 
angle between the segment $PP^*$ and a horizontal line.
The difference between the
old cost and the new cost can be expressed as
\bel{onc1} |P - P^*| \left(\sigma + \sum_{j=1}^{n}\kappa_j \right)^\alpha + \sum_{j=1}^{n}\kappa_j^\alpha |P^* - P_j|
- \sigma^\alpha \cos{\theta_a} |P-P^*| - \sum_{j=1}^{n}\left( \sum_{i=1}^{j}\kappa_i \right)^\alpha
|P_{j+1} -P_j|, \eeq
where, for notational convenience,  we set $P_{n+1} \doteq P$.
According to (\ref{onc1}) we can write
\bel{ABn} \hbox{\rm [old cost]} - \hbox{\rm [new cost]}~=~A+ S_n\,,\eeq
where
\bel{Adef} A~\doteq~|P-P^*| \left[ \left(\sigma + \sum_{j=1}^{n}\kappa_j \right)^\alpha - \sigma^\alpha \cos{\theta_a} \right] +
\left( \sum_{j=1}^{n}\kappa_j \right)^\alpha \Big( |P^* - P_1| - |P - P_1| \Big) ,\eeq
\bel{Sn}
S_n~=
~\sum_{j=1}^{n}\kappa_j^\alpha |P^* - P_j| -
 \Big( \sum_{j=1}^{n}\kappa_j \Big)^\alpha
\Big( |P^* - P_1| - |P_{n+1} - P_1| \Big) - \sum_{j=1}^{n}\Big( \sum_{i=1}^{j}\kappa_i \Big)^\alpha
|P_{j+1} -P_j|.
 \eeq
\v
{\bf 2.} 
Notice that the quantity $A$ in (\ref{Adef}) would 
describe the difference  in the costs
if all the mass $\kappa=\kappa_1+\cdots+\kappa_n$
were flowing through the point $P_1$. Using  Lemma \ref{l:41},
we can thus choose $P=P_1$ close enough to $P^*$ 
such that this difference is strictly positive.
More precisely, for a fixed $\kappa>0$, we claim that one can achieve the lower bound
\bel{lb1}\bega{rl}A&\ds \geq~ |P-P^*| \left[(\sigma + \kappa)^\alpha - \sigma^\alpha \cos{\theta_a} + \kappa^\alpha
\cos(\theta_a + \theta_1) - {\kappa^\alpha \over 2} {|P-P^*| \over |P_1 - P^*|}\, \right]
\\[3mm]
&\geq~|P_1-P^*|\cdot\delta(\theta_1,\kappa)~>~0.\enda
\eeq
Indeed,  the last two terms within the square brackets in (\ref{lb1}) are derived  from
$$ \bega{rl} |P^*-P_1| - |P-P_1| &=~\ds |P^*-P_1|\left[1 - \sqrt{1 - 2{|P-P^*| \over |P^*-P_1|}\cos(\theta_a+\theta_1)
+ {|P-P^*|^2 \over |P^*-P_1|^2} } \right] \\
&\geq~\ds |P^*-P_1|\left[1 - \left(1 - {|P-P^*| \over |P^*-P_1|}\cos(\theta_a+\theta_1)
	+ {|P-P^*|^2 \over 2|P^*-P_1|^2} \right) \right]. \enda $$
Moreover, since we have the strict inequalities
\bel{t1s}\left\{\bega{rll}&\theta_1\,<\,{\pi\over 2} \qquad& \hbox{if}\quad \alpha\geq {1\over 2}\,,\\[4mm]
 &\theta_1<{\pi\over 2}-\theta_0\qquad &\hbox{if}\quad \alpha< {1\over 2}\,,\enda
 \right.\eeq
the same argument used the proof of (\ref{c3})  in Lemma~\ref{l:41} now yields the
strict inequality
\bel{c33} \cos{\theta_1} ~< ~{ \lambda^{2\alpha} + (1-\lambda)^{2\alpha} -1 \over 2 \lambda^\alpha (1-\lambda)^\alpha }\,. \eeq
Given $\kappa>0$ and $P_1$,  
we can then choose $P$ close enough to $P^*$ 
so that 
\begi
\item  the terms within the square brackets in (\ref{lb1}) is strictly positive,
\item the ratio $|P-P^*|/|P_1-P^*|$ is small but uniformly positive, as long as
$\theta_1$ remains bounded away from ${\pi\over 2}$ or from ${\pi\over 2}-\theta_0$
respectively, in the two cases considered in (\ref{t1s}).
\endi
This proves our claim (\ref{lb1}).
\v
{\bf 3.} To complete the proof of the lemma, it remains to prove that  $S_n \ge 0$.
This will be proved by induction on $n$.   Starting from (\ref{Sn}) and using the inequalities 
$$ |P_n - P_1|~ \leq~ |P^* - P_1|,\qquad\quad \Big( \sum_{i=1}^n\kappa_i \Big)^\alpha~\leq~\kappa_n^\alpha + \Big( \sum_{i=1}^{n-1}\kappa_i \Big)^\alpha,$$
 we obtain
\bel{Sn1}
\bega{rl} S_n&=\ds~\sum_{j=1}^{n}\kappa_j^\alpha |P^* - P_j|  -\Big( \sum_{j=1}^{n}\kappa_j \Big)^\alpha
\underbrace{\bigl( |P^* - P_1| - |P_n - P_1| \bigr) }_{\ge 0}
- \sum_{j=1}^{n-1}\Big( \sum_{i=1}^{j}\kappa_i \Big)^\alpha |P_{j+1} -P_j|\\[4mm]
&\ds\geq ~\sum_{j=1}^{n-1}\kappa_j^\alpha |P^* - P_j|  -\Big( \sum_{j=1}^{n-1}\kappa_j \Big)^\alpha\bigl( |P^* - P_1| - |P_{n-1} - P_1| \bigr)
- \sum_{j=1}^{n-2}\Big( \sum_{i=1}^{j}\kappa_i \Big)^\alpha |P_{j+1} -P_j|
\\[4mm]
&\qquad\ds +\kappa_n^\alpha  |P^*-P_n| - \kappa_n^\alpha 
\Big( |P^* - P_1|- |P_n - P_1|  \Big)+\Big( \sum_{i=1}^{n-1}\kappa_i \Big)^\alpha |P_{n} -P_{n-1}|
\\[4mm]
&= ~\ds S_{n-1} +  \kappa_n^\alpha\Big(  |P^*-P_n| -|P^* - P_1|+ |P_n - P_1| \Big)
~\geq~S_{n-1}\,.
\enda \eeq
Repeating this same argument, by induction we obtain
$$S_n~\geq~S_{n-1}~\geq~\cdots~\geq S_1\,.$$
Observing that
$$ S_1 ~=~  \kappa_1^\alpha |P^* - P_1| - \kappa_1^\alpha \Big( |P^* - P_1| - |P_2 - P_1| \Big)- \kappa_1^\alpha
|P_2 - P_1|~ =~ 0, $$
we complete the proof of the lemma. \endproof

\begin{figure}[ht]
	\centerline{\hbox{\includegraphics[width=10cm]{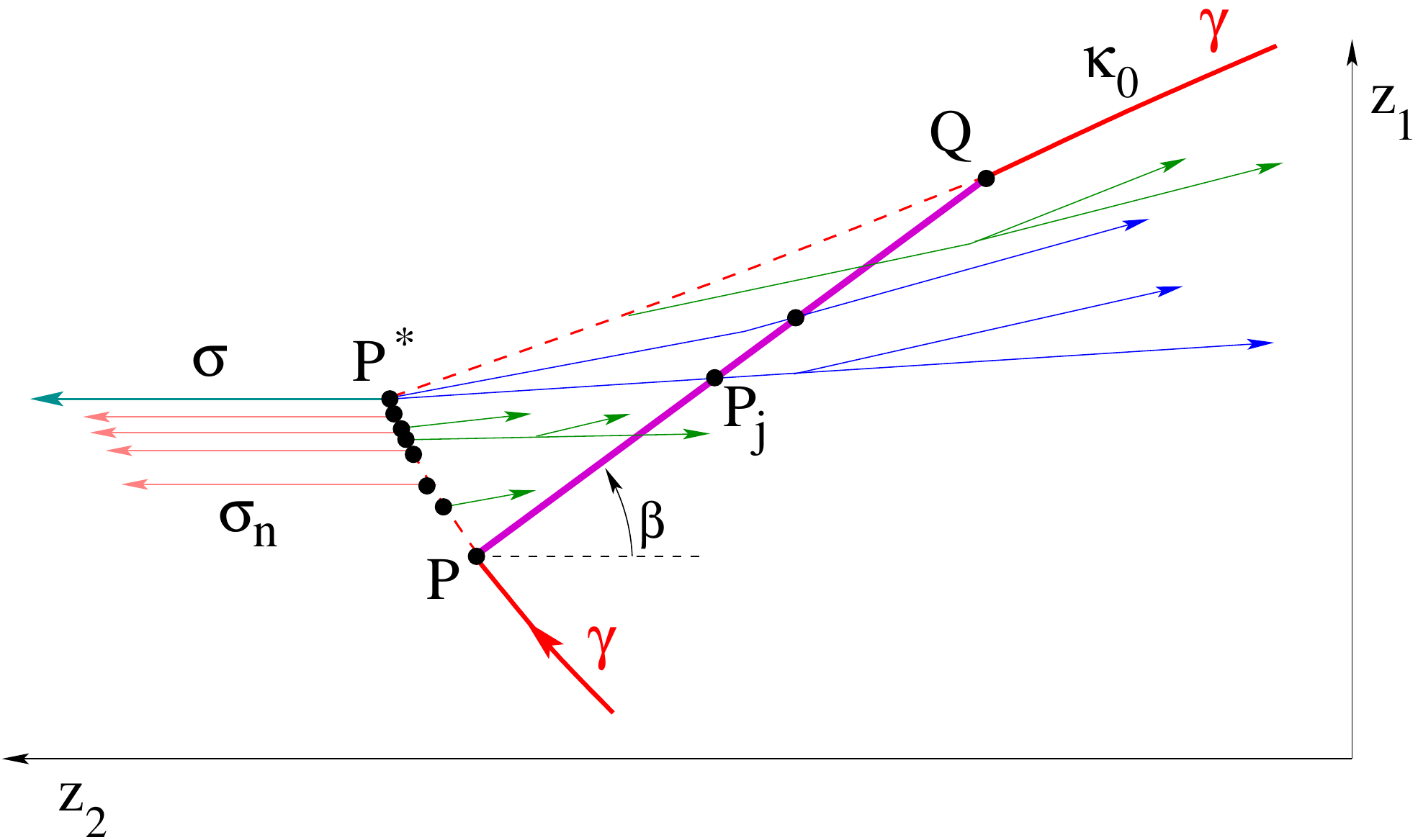}}}
	\caption{\small In the fully general situation, we have additional branches bifurcating
	to the left of $\gamma$ between $P$ and $P^*$, and to the right of $\gamma$ 
	at any point between $P$ and $Q$.   In addition, there can be an additional absolutely
	 continuous source along the arc $PP^*$.} 
	\label{f:ir131}
\end{figure}

\v
We now consider the most general situation, shown in 
Fig.~\ref{f:ir131}.  Differently from the setting of Lemma~\ref{l:42}, various 
scenarios must be considered.
\begi
\item In addition to the horizontal path bifurcating to the left of $P^*$ with flux $\sigma$, there can be
countably many additional horizontal branches bifurcating 
to the left of $\gamma$, below $P^*$.
We shall denote by $\sigma_n$, $n\geq 1$, the fluxes through these branches, at the bifurcation points.
\item  There can be countably many distinct branches bifurcating 
to the right of $P^*$, say with fluxes $\kappa^*_j$, $j\geq 1$.
\item Furthermore, there can be countably many additional branches
bifurcating to the right of $\gamma$, at points close to $P^*$.
We shall denote by $\kappa'_i$, $i\geq 1$, the fluxes through these branches, at the bifurcation points.
\item Finally, the measure $\mu$ could concentrate a positive mass 
along the arc $PP^*$.
\endi
We observe that, by optimality, all particle trajectories  to the right of $\gamma$
move in the right-upward direction. Namely, setting
$\chi(\xi,t)= (z_1(\xi,t), \, z_2(\xi,t))$, for these paths we have
$$\dot z_1(\xi,t)\,\geq\,0,\qquad \dot z_2(\xi,t)\,\leq\,0.$$

We now construct a ``bypass", 
choosing a segment $PQ$ with endpoints
both lying on the  curve $\gamma$, making an angle $\beta$ with the horizontal 
direction such that
\bel{bb*}\beta^*~<~\beta~<~{\pi\over 2} - \theta_0\,.\eeq
Here $\beta^*$  denotes the angle between the segment $P^*Q^*$ and a horizontal line.

Given $\ve>0$, 
we can choose $N\geq 1$ large enough so that, among the 
branches bifurcating from $P^*$, one has
\bel{sm1} \sum_{j>N} \kappa^*_j~<~\ve.\eeq
Moreover, by choosing $Q$ sufficiently close to $P^*$, the following can be achieved:

\begi
\item[(i)] The total flux along the horizontal branches bifurcating 
to the left of $\gamma$ below $P^*$ satisfies
\bel{sm2}\sum_{n\geq 1} \sigma_n~<~\ve.\eeq
\item[(ii)] The total flux along the branches bifurcating 
to the right of $\gamma$ between $P$ and $P^*$, and between 
$P^*$ and $Q$ satisfies
\bel{sm3}\sum_{i\geq 1} \kappa'_i~<~\ve.\eeq
\item[(iii)] For each $j=1,\ldots,N$, there exists a path $\gamma_j$ connecting
$P^*$ with a point $P_j$ on the segment $PQ$, along which
the flux remains $\geq \kappa_j~\geq ~\kappa_j^*-(\ve/N)$.
Here we denote by $\kappa_j$ the flux reaching  $P_j$.

In other words, even if the $j$-th branch through $P^*$ further bifurcates,
most of the particles along this branch cross the segment $PQ$ at the same 
point $P_j$.
\item[(iv)] The total mass of $\mu$ along $\gamma$, between 
$P$ and $P^*$ is $<\ve$.
\endi

We estimate the difference in the new cost produced by these additional 
branches.
Call $P=(p_1, p_2)$, $Q=(q_1, q_2)$.

\begi
\item The additional mass on the left branches, together with the mass
of $\mu$ present between $P$ and $P^*$ now travels along
a horizontal line through $P$. By (i) and (iv)  this mass is $<2\ve$.
Hence:
\bel{e11}\hbox{[additional cost]} ~\leq~(2\ve)^{1-\alpha} (z_2^*- z_2).\eeq
\item
The additional mass bifurcating to the right of $\gamma$, not crossing
the segment $PQ$ at one of the finitely many points 
$P_1,\ldots, P_N$ is $< 3\ve$.   The additional 
cost in transporting this mass from $P$ to some point  between $P$ and $Q$
satisfies
\bel{e22}\hbox{[additional cost]} ~< \kappa_0^{\alpha-1} \cdot 3\ve |P-Q|.
\eeq
\endi

We now use Lemma~\ref{l:42}.  Combining (\ref{saving}) with
(\ref{e11})-(\ref{e22}) we obtain
\bel{sav}
\hbox{\rm [old cost]} - \hbox{\rm [new cost]}~\geq~|P_1-P^*|\cdot \delta(\theta_1,\kappa) -(2\ve)^{1-\alpha} |P-P^*| -\kappa_0^{\alpha-1} \cdot 3\ve |P-Q|. \eeq
By choosing $\ve>0$ small enough, the right hand side of (\ref{sav}) is strictly positive.   Hence the configuration with $P^*\not= 0$ is not optimal.  This completes the proof
of Theorem~\ref{t:1}.
\endproof

\section{The case $d=2$, $\alpha=0$}
\label{s:6}
\setcounter{equation}{0}

We give here a proof of Theorem~\ref{t:2}.
\v
{\bf 1.}  Assume that there exists a unit vector $\bfw^*\in \R^2$ such that
$$K~=~~ \int_{\bfn\in S^1} \Big|
\langle \bfw^*, \bfn\rangle\Big|\, \eta(\bfn)\, d\bfn~>~c.$$
Let $\bfv= (\cos \beta, \sin\beta)$ be a unit vector perpendicular to $\bfw^*$,
with $\beta\in[0,\pi]$.
Let $\mu$ be the measure supported on the segment
$\{r\bfv\,;~r\in [0,\ell]\}$, with constant density $\lambda$ w.r.t.~1-dimensional Lebesgue
measure.

Then the payoff achieved by $\mu$ is estimated by
\bel{po8}\bega{rl}\ds\S^\eta(\mu)-c\I^0(\mu)&=~\ds
\ell\cdot \int_{S^1} \left( 1- \exp\Big\{ - {\lambda\over \bigl|\langle \bfw^*, \bfn\rangle
\bigr|}\Big\}\right)\, \Big|
\langle \bfw^*, \bfn\rangle\Big|\, \eta(\bfn)\, d\bfn - c\,\ell\\[4mm]
&\ds\geq~
\ell\cdot (1-e^{-\lambda})\, \int_{S^1}\Big| \langle \bfw^*, \bfn\rangle
\Big|\, \eta(\bfn)\, d\bfn - c\,\ell
\\[4mm]
&=~\Big[(1-e^{-\lambda})\,K- c\Big]\,\ell.\enda\eeq
By choosing $\lambda>0$ large enough, the first factor on the right hand side of 
(\ref{po8}) is strictly positive.   Hence, by increasing the length $\ell$, we can render
the payoff arbitrarily large.
\v 
{\bf 2.} Next, assume that $K\leq c$. 
Consider any Lipschitz curve $s\mapsto \gamma(s)$, parameterized by 
arc-length $s\in [0,\ell]$.   Then, for any measure $\mu$ supported on $\gamma$,
 the total amount of sunlight  from the direction $\bfn$ captured by $\mu$  
 satisfies the estimate
$$\S^\bfn(\mu)~\leq~\int_0^\ell 
\Big|\langle \dot\gamma(s)^\perp ,\, \bfn\rangle\Big|\, ds.$$
Indeed, it is bounded by the length of the projection of $\gamma$ on the 
line $E_\bfn^\perp$ perpendicular to $\bfn$.
Integrating over the various sunlight directions, one obtains
$$\S^\eta(\mu)~\leq~\int_0^\ell \int_{S^1}
\Big|\langle \dot\gamma(s)^\perp ,\, \bfn\rangle\Big|\,\eta(\bfn)\, d\bfn\, ds~\leq~K\,\ell.$$

More generally, $\mu= \sum_i \mu_i$ can be the sum of countably many 
measures supported on Lipschitz curves $\gamma_i$.
In this case, since the sunlight functional is sub-additive, one has
$$\S^\eta(\mu)~\leq~\sum_i \S^\eta(\mu_i)~\leq~\sum_i K \ell_i\,.$$
Hence 
$$\S^\eta(\mu)-c\I^0(\mu)~\leq~\sum_i K\ell_i - c \sum_i\ell_i~\leq~0.$$
This concludes the proof of case (ii) in Theorem~\ref{t:2}.
\endproof

\v

{\bf Acknowledgments.}
The research of the first author was partially supported by NSF with  
grant DMS-1714237, ``Models of controlled biological growth".
The research of the second author was partially supported by a grant from the
U.S.-Norway Fulbright Foundation.

\v

\end{document}